\documentclass[11pt]{amsart}
\usepackage[titletoc,title]{appendix}
\usepackage{amsmath,amsthm,amsfonts,amssymb,mathrsfs,epsfig,bm,bbm}

\usepackage[titletoc,title]{appendix}
\usepackage{amsmath,amsthm,amsfonts,amssymb,mathrsfs,epsfig,bm}

\usepackage[usenames]{color}

\usepackage[colorlinks=true, urlcolor={black}]{hyperref}
\usepackage{dsfont}
\usepackage{cite}

\newtheorem{theorem}{Theorem}
\newtheorem{lemma}{Lemma}

\theoremstyle{remark}

\newtheorem{remark}{Remark}

\def\Z{\mathbb{Z}}

\def\R{\mathbb{R}}

\def\P{\mathbb{P}}
\def\E{\mathbb{E}}

\renewcommand{\phi}{\varphi}
\renewcommand{\epsilon}{\varepsilon}
\newcommand{\1}{{\text{\Large $\mathfrak 1$}}}

\definecolor{mygray}{gray}{0.9}

\newcommand{\keyword}[1]{ \noindent {\footnotesize
             {\small \em Keywords and phrases.} {\sc #1} } }
\newcommand{\ams}[2]{  \noindent {\footnotesize
             {\small \em AMS {\rm 2000} subject classifications.
             {\rm  {\sc #1};  {\sc #2}} } } }
\newcommand{\acknowledgments}{\par{\noindent\bf Acknowledgements. }}

\begin{document}

\begin{abstract}
The greedy walk is a walk on a point process that always  moves 
from its current position to the nearest not yet visited point.
We consider here various  point processes on two lines.
We look first at the greedy walk on two independent one-dimensional Poisson 
processes placed on two intersecting lines and prove that the greedy walk
 almost surely does not visit all points.
When a point process is defined on two parallel lines, the result depends on the definition of the process:
If each line has a copy of the same realisation of a homogeneous Poisson point process, 
then the walk almost surely does not visit all points of the process.
However, if each point of this process 
is removed with probability $p$ from either of the two lines, independently of the other points,
then
the walk almost surely visits all points.
Moreover, the greedy walk on two parallel lines, where each line has a copy of the same realisation of a homogeneous Poisson point process,  
but one copy is shifted by some small
$s$, almost surely  visits all points.

\vspace*{2mm}

\keyword{Poisson point process; greedy walk}

\ams{60K37}{60G55,60K25}
\end{abstract}

\title{Greedy walks on two lines}

\author{\sc Katja Gabrysch}

\maketitle


\section{Introduction}

Consider a simple point process $\Pi$ in a metric space $(E,d)$.
We think of $\Pi$ as a collection of points (the support of the measure) 
and 
we use the notation $\vert \Pi\cap B\vert$
 to indicate the number of points 
on the Borel set $B \subset E$.
If $x \in E$, the notation $x \in \Pi$ is used instead of $\Pi\cap \{x\}\neq \emptyset$.

We define a greedy walk on $\Pi$
as follows. 
The walk 
starts from some point $S_0 \in E$ and always moves
on the points of $\Pi$ by picking the point closest to its current
position that has not been visited before. 
Thus a sequence $(S_n)_{n\geq 0}$ is defined recursively by
\begin{equation}
\label{defwalk}
S_{n+1} = \arg\min \big\{d(X, S_n ) :X \in \Pi, ~ X\notin \{S_0, S_1, \ldots, S_n\} \big\}.
\end{equation}

The greedy walk is a model in queuing systems where the points of the process represent
positions of customers and the walk represents a server moving towards customers.
Applications of such a system can be found, for example, in telecommunications and computer networks or 
transportation.
As described in \cite{BFL11}, the model of a greedy walk on a point process can be defined in various ways and on different spaces.
For example, Coffman and Gilbert \cite{CG87} and Leskel\"{a} and Unger \cite{LU12} study
a dynamic version of the greedy walk on a circle with new customers arriving to the system 
according to a Poisson process.

The greedy walk defined  as in \eqref{defwalk}
on a homogeneous Poisson process on $\R$ 
almost surely does not visit all points.
More precisely, the expected number of times the walk jumps over $0$ is $1/2$ \cite{Katja15}.
Foss {\em et al.} \cite{FRS} and Rolla {\em et al.} \cite{RST14} study two modifications of this model where
they introduce some extra points on the line, which they call ``rain'' and ``dust'', respectively. 
Foss {\em et al.} \cite{FRS} consider a space-time model,
starting with a Poisson process at time 0. The positions and times of arrival of new points
are given by a Poisson process on the half-plane.
Moreover, the expected time that the walk spends at a point is 1.
In this case the walk, almost surely, jumps over the starting point finitely many times
and the position of the walk diverges logarithmically in time.
Rolla {\em et al.} \cite{RST14} assign to the points of a Poisson process one or two marks at random. 
The walk always moves to the point closest to the current position which still has at least one mark left and then
removes exactly one mark from that point.
The authors show that introducing points with two marks will force the walk to change
sides infinitely many times.
Thus, unlike the walk on a Poisson process with single marks, the walk here almost surely visits all points of the point process.

There is not much known about the behaviour of the greedy walk on a homogeneous Poisson process in higher dimensions.
For example, it is an open problem whether 
the greedy walk on the points of a homogeneous Poisson process on $\R^2$ visits all points \cite{RST14}.

In this paper, we study various point processes defined on the union of two lines $E\subset\R^2$,
where the distance function $d$ on $E$ is the Euclidean distance,
$d\left((x_1,y_1),(x_2,y_2)\right)=\sqrt{(x_1-x_2)^2+(y_1-y_2)^2}$.
For all point processes $\Pi$ considered in this paper, every step of the greedy walk on $\Pi$ is almost surely 
uniquely defined, that is for every $n\geq 0$, there is, almost surely, only one point for which the minimum \eqref{defwalk} is obtained.

We study first a point process $\Pi$ on two lines intersecting at $(0,0)$
with independent homogeneous Poisson processes on each line.
The greedy walk starts from $(0,0)$.
When the walk visits a point that is far away from $(0,0)$, then 
the distance to $(0,0)$ and to any point on the other line is large.
Thus, the probability of changing lines or crossing $(0,0)$ is small.
In Section \ref{intersecting} we show by using the Borel-Cantelli lemma that almost surely
the walk crosses $(0,0)$ or changes  lines only finitely many times,
which implies that almost surely the walk does not visit all points of $\Pi$.

Thereafter, we look at the greedy walk on two parallel lines
at a fixed distance $r$, $\R\times\{0,r\}$, with a point process on each line.
The behaviour here depends on the definition of the process.
The first case we study is a process $\Pi$ 
consisting of two identical copies of 
a homogeneous Poisson process on $\R$, 
that are placed 
on the parallel lines. 
We show in Section \ref{parallel_same_process} that the greedy walk does not visit all the points of $\Pi$, 
but it visits all the points on one side of the vertical line $\{0\}\times \R$ and just finitely many points on the other side.

In the second case, we modify the definition of the process 
above by deleting exactly 
one of the copies of each point 
with probability $p>0$, independently from the other points,
and the line is chosen with probability $1/2$.  
In particular, if $p=1$ we have two independent Poisson processes on these lines.
For any $p>0$, the greedy walk almost surely visits all points.
The reason is that 
the greedy walk skips some of the points when it goes away from the vertical line $\{0\}\times \R$ 
and those points will force the walk to return and cross the vertical line 
$\{0\}\times \R$ infinitely many times. 
We prove this in Section \ref{par_indep} using arguments from \cite{RST14}.

The greedy walk also visits  all points of $\Pi$ in the case when $\Pi$ consists of two identical copies of 
 a homogeneous Poisson process on $\R$ 
where one copy is shifted by $\vert s\vert<r/\sqrt{3}$.
This is discussed in Section \ref{shift}. 
Note that all results are independent of the choice of $r$.  

For the greedy walk on a homogeneous Poisson process on $\R$, 
with single or double marks assigned to the points,
Rolla {\em et al.} \cite{RST14} show that, even though the walk visits all the points, 
the expected first crossing time of $0$ is infinite.
One can in a similar way show analogous results for the greedy walk on the processes on two parallel lines defined in Sections \ref{par_indep} and \ref{shift},
but that is not included in this paper.

\section{Two intersecting lines}
\label{intersecting}

Let $E=\{(x,y)\in\R^2:y=m_1x \text{ or } y=m_2x\}$ where $m_1,m_2\in \R$, $m_1\neq m_2$.
The point $(0,0)$ divides $E$ into four half-lines.
Let $d$ be the Euclidean distance and
denote the distance of a point $(x,y)$ from the origin by $| |(x,y)| |=d((0,0),(x,y))$.

Let $\Pi^1,\Pi^2$ be two independent Poisson processes on $\R$ with rate 1. 
Then, for any $a,b\in \R$ let
$$\Big\lvert \Pi\cap\Big\{(x,m_ix):x\sqrt{1+m_i^2}\in (a,b)\Big\}\Big\rvert=\big\lvert\Pi^i\cap(a,b)\big\rvert,$$
so that the distances between the points of $\Pi^1$ and $\Pi^2$ are preserved in $E$.
The greedy walk $(S_n)_{n\geq 0}$ on $\Pi$ defined by \eqref{defwalk} starts from $(0,0)$.

\begin{theorem}
\label{t_three}
 Almost surely, the greedy walk does not visit all points. 
More precisely, the greedy walk almost surely visits only finitely many points on three half-lines.

\end{theorem}

\proof
Let $B_{(0,0)}(R)$ be a ball in $\R^2$ of radius $R$ around point $(0,0)$. Then
$$\vert\Pi\cap B_{(0,0)}(R)\vert= \vert \Pi^1\cap (-R,R)\vert+\vert \Pi^2 \cap (-R,R)\vert<\infty \quad \text{ a.s.}$$
Thus, $\lim_{n \to \infty}| |S_n|| = \infty$ almost surely.
To show that the walk does not visit all points of $\Pi$, 
it suffices to prove that 
the sequence $(S_n)_{n\geq 0}$ changes half-lines only finitely many times. 

We can define a subsequence of the times when the walk changes half-lines 
as follows:
\begin{align*}
 j_0&=0, \\
 j_n&= \inf\left\{k>j_{n-1}: S_{k} \text{ and } S_{k+1} 
\text{ are not on the same half-line }\right.\\ 
&\qquad\qquad\qquad\qquad\text{and } 
 \left. | |S_{k}||>\max\{n, ||S_{j_{n-1}}||\}\right\},\\
k_0&=0,\\
k_n&=\inf\left\{k\leq j_n:S_{k},S_{k+1},\dots,S_{j_n} 
\text{ are on the same half-line}\right\},
\end{align*}
where $\inf \emptyset=\infty$.
Let $(U_{n})_{n\geq 0}$ and $(V_{n})_{n\geq 0}$  be the corresponding subsequences of $(S_n)_{n\geq 0}$, that is,
$$U_n=S_{k_n}\text{ and } V_n=S_{j_n},$$
when $j_n,k_n<\infty$, and $U_n=(1,1)$, $V_n=(0,0)$ otherwise.
Moreover, define the events
$$B_n=\{||U_{n}||\leq||V_{n}||\} \text{ and }C_n=\{||U_{n}||> ||V_{n}||\}.$$

If the greedy walk changes half-lines infinitely often, then $j_n<\infty$ for all $n$
and exactly one of the events $B_n$ and $C_n$ occurs for each $n$. 

Assume that 
$C_n$ occurs for some $n$ such that $j_n<\infty$.
Let $X=S_{k_n-1}$ be the last visited point before $U_n$. 
Then, 
by the definition of the sequence $(V_n)_{n\geq 0}$,
$\vert |X|\vert\leq \max\{n,\vert |V_{n-1}|\vert\}<\vert |V_{n}|\vert$.
From $||U_n||>||V_n||>||X||$ it follows that 
$d(X,V_n)<d(X,U_n)$, which contradicts the definition of the greedy walk, and therefore 
also the assumption that $C_n$ occurs.

Assume now that $B_n$ occurs for some $n$ such that $j_n<\infty$.
Denote by $\alpha$ the acute angle between the lines $y=m_1x$ and $y=m_2x$.
By the definition of the sequence $(V_n)_{n\geq 0}$, 
the walk up to time $j_n$ never changed lines from a point whose distance from the origin is greater than $||V_{n}||$. 
Moreover, because of the assumption $||U_{n}||\leq ||V_{n}||$, between times $k_n$ and $j_n$ the walk never visited
the points further away than $V_{n}$ on the corresponding half-line.
Thus, the points on this half-line 
outside $B_{(0,0)}(||V_{n}||)$
are not yet visited.
The distance from $V_{n}$ to the other line is $||V_{n}||\sin \alpha$.
Since the walk changes lines after visiting the point $V_{n}$, 
we can conclude that there are no unvisited points of $\Pi$ in $B_{V_n}(||V_{n}||\sin \alpha)$ at time $n$.
Hence, for $n\geq n_0$, where $n_0$ is such that $n_0\sin \alpha>1$, we have
\begin{align*}
B_{n} 
&\subset
\bigcup_{i=1}^{2}\left\{\Pi^{i}\cap \left( ||V_{n}||,||V_{n}||(1+\sin \alpha)\right)=\emptyset\right\}
\cup\left\{\Pi^{i} \cap \left( -||V_{n}||(1+\sin \alpha),-||V_{n}|| \right)=\emptyset\right\}\\
&\subset 
\bigcup_{i=1}^{2}\bigcup_{R=n}^{\infty}\left\{\Pi^{i} \cap \left( R+1,R(1+\sin \alpha) \right)=\emptyset\right\}
\cup\left\{\Pi^{i} \cap \left( -R(1+\sin \alpha),-(R+1) \right)=\emptyset\right\}.
\end{align*}
Then,
$$ \P(B_{n})\leq 4\sum_{R=n}^{\infty}e^{-R\sin \alpha+1}=\frac{4e^{-n\sin \alpha+1}}{1-e^{-\sin \alpha}} $$
and 
$$\sum_{n=n_0}^\infty  \P(B_{n})\leq \frac{4e^{-n_0\sin \alpha+1}}{(1- e^{-\sin \alpha})^2}<\infty.$$
Hence, by the Borel-Cantelli lemma,  $ \P\left(B_{n} \text{ for infinitely many }n\geq 1 \right)=0$.

Since the event $C_n$ does not occur for any $n$ such that $j_n<\infty$
and $B_n$ occurs almost surely finitely many times,  
there exists, almost surely,  $n_0$ such that $j_n=\infty$ for $n\geq n_0$.
Thus, the greedy walk changes lines finitely many times.  
\qed

\begin{remark}
The theorem holds also when $\Pi^1$ and $\Pi^2$ are Poisson processes with different rates.
Moreover, we can generalise this theorem as follows.
 Let $E$ be a space of finitely many intersecting lines 
(every two lines are intersecting, but not necessarily all in the same point)
and $\Pi$ is a point process on $E$ consisting of a homogeneous Poisson process (with possibly different rates) 
on every line.
Then the greedy walk starting from a point in $E$ does not visit all points of $\Pi$.
Similarly as above, one can show that the walk changes lines finitely many times, and therefore, visits 
finitely many points on all but one line.
\end{remark}

\section{Two parallel lines with the same Poisson process}
\label{parallel_same_process}
Let $E = \R\times \{0,r\}$ and let $d$ be the
Euclidean distance.
Sometimes we refer to the lines $\R\times\{0\}$ and $\R\times\{r\}$ as line $0$ and line $r$, respectively.

We define a point process $\Pi$ on $E$ in the following way.
Let $\widehat \Pi$ be a homogeneous Poisson process on $\R$ with rate $1$ 
and let $\vert \Pi\cap (B\times\{0\})\vert=\vert\Pi\cap (B\times\{r\})\vert=\vert\widehat \Pi\cap B\vert$ for all Borel sets $B\subset\R$.
Denote the points of the process $\widehat \Pi$ by
$$\ldots<X_{-2}<X_{-1}\leq 0<X_1<X_2<\dots.$$
The greedy walk on $\Pi$ starts from the point $S_0=(0,0)$, which is with probability $1$ not a point of $\Pi$.

\begin{figure}[h]
    \centering
    \includegraphics[width=\textwidth]{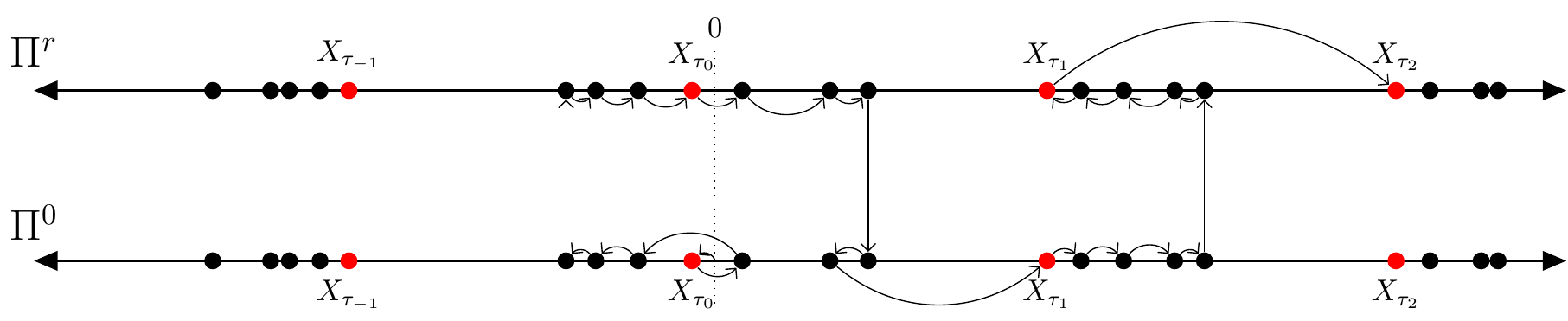}
    \caption{ 
The points closest to $0$ in each cluster form the sequence $(X_{\tau_i})_{i\in \Z}$.
Once the greedy walk enters a cluster, it visits successively all the points of that cluster
before it moves to the next cluster.  
}
    \label{clusters}
\end{figure}

\begin{theorem}
\label{two_copies}
Almost surely, the greedy walk does not visit all points of $\Pi$.
More precisely, 
the greedy walk almost surely visits all points on one side of the vertical line $\{0\}\times \R$
and finitely many points on the other side.
\end{theorem}
\proof
We can divide the points $(X_i)_{i\in\Z}$ into clusters, 
so that successive points in the same cluster have a distance less than $r$ 
and the distance between any two points in different clusters is greater than $r$.
Let $(\tau_i)_{i\in \Z}$ be the indices of the closest point to 0 in each cluster. 
More specifically,  $\tau_0=-1$ if $\vert X_{-1}\vert <X_1$ and $\tau_0=1$  otherwise,
$(\tau_i)_{i\geq 1}$ is the unique sequence of integers such that $X_{\tau_i}-X_{\tau_i-1}>r$
and $X_{k}-X_{k-1}\leq r$ for $\tau_{i-1}<k<\tau_i$, and, similarly, $(\tau_{-i})_{i\geq 1}$
is a sequence of integers such that $X_{\tau_{-i}+1}-X_{\tau_{-i}}>r$ and $X_{k+1}-X_{k}\leq r$ for $\tau_{-i-1}<k<\tau_{-i}$.
Moreover, we call the cluster containing the point $X_{\tau_i}$ cluster $i$.
See Figure \ref{clusters} for an example of clusters $-1,0,1$ and $2$. 
The points $\{X_{\tau_{-1}},X_{\tau_{0}}, X_{\tau_{1}}, X_{\tau_{2}}\}\times\{0,r\}$  are marked with gray colour.

The greedy walk starting from $(0,0)$ 
visits several points around $(0,0)$ and then moves to line $r$ from one of the outermost points of cluster $0$ 
on line $0$.
Then it visits all points of cluster $0$ on line $r$ and if there are points left it changes lines again to visit the remaining points
on line $0$ before moving to the next cluster.
Later, when the greedy walk is in cluster $i$, $i\neq 0$, 
it always visits first a point at $X_{\tau_i}$, that is $(X_{\tau_i},0)$ or $(X_{\tau_i},r)$. 
Then the walk visits successively all the points of cluster $i$ on the same line and it changes lines
at the other outermost point of the cluster.
Thereafter the walk visits the corresponding points on the other line
 in reverse order until it reaches a point at $X_{\tau_i}$.
Thus, the walk visits all points of  cluster $i$ consecutively and it ends at the starting position $X_{\tau_i}$, 
but on the other line.
Therefore, to know whether the walk visits all points, it is enough to know 
the positions of the points in cluster  $0$
and the position of the points at $X_{\tau_i}, i\neq 0$.
Since the points of cluster $0$ almost surely do not change the asymptotic behaviour of the walk, 
for the proof we look at the greedy walk $(\widetilde S_n)_{n\geq 0}$ on 
$(X_{\tau_i})_{i\in\Z}$ with $0$ as starting point.
Note that the distances $X_{\tau_{i+1}}-X_{\tau_{i}}$, $i\in\Z\setminus \{-1,0\}$,
are independent and identically distributed random variables with finite expectation.

To visit all points $(X_{\tau_i})_{i\in\Z}$, the greedy walk needs to cross over $0$ infinitely many times.
Let $A_{m}$ be the event that the walk $(\widetilde S_n)_{n\geq 0}$ crosses $0$ after visiting a point in the interval $(rm,r(m+1))$, 
that is
$$A_{m} =\{\exists\ n : rm \leq \widetilde S_n < r(m+1) \text{ and } \widetilde S_{n+1}< 0\}.$$
This can be written as
\begin{align*}
 A_{m}&=\{\exists\ n,i:rm \leq X_{\tau_i}< r(m+1),\  \widetilde S_n=X_{\tau_i} \text{ and } \widetilde S_{n+1}<  0 \}\\
&\subset\{\exists\ i: rm \leq X_{\tau_i}< r(m+1),\ X_{\tau_{i+1}}-X_{\tau_i}>rm\}\\
&=\bigcup_{i\geq 0} \{rm \leq X_{\tau_i}< r(m+1),\ X_{\tau_{i+1}}-X_{\tau_i}>rm\}.
\end{align*}
For $i\geq 1$, the random variable $X_{\tau_{i+1}}-X_{\tau_i}$ is independent of $X_{\tau_i}$ and 
it has the same distribution as $X_{\tau_{2}}-X_{\tau_1}$.
Moreover, $X_{\tau_i}\in [rm,r(m+1))$ for at most one $i$.
Thus,  we have 
\begin{align*}
 \P(A_{m})&\leq \P(rm \leq X_{\tau_0}< r(m+1))+\sum_{i=1}^{\infty} \P(rm \leq X_{\tau_i}< r(m+1),\ X_{\tau_{i+1}}-X_{\tau_i}>rm)\\
&\leq \P(rm \leq X_{\tau_0}< r(m+1))+\sum_{i=1}^{\infty} \P(rm \leq X_{\tau_i}< r(m+1))\P(X_{\tau_{2}}-X_{\tau_1}>rm)\\
&= \P(rm \leq X_{\tau_0}< r(m+1))+\P(X_{\tau_{2}}-X_{\tau_1}>rm) \sum_{i=1}^{\infty} \E \left(\1_{\left\{rm 
\leq X_{\tau_i}< r(m+1)\right\}}\right)\\ 
 &\leq \frac{1}{2}e^{-2rm}(1-e^{-2r})+\P\left(X_{\tau_{2}}-X_{\tau_1}>rm\right).
\end{align*}
Then 
\begin{align*}
 \sum_{m=1}^{\infty} \P(A_{m})&\leq  \sum_{m=1}^{\infty}\frac{1}{2}e^{-2rm}(1-e^{-2r})+ \sum_{m=1}^{\infty}\P(X_{\tau_{2}}-X_{\tau_1}>rm)\\
&<\frac{1}{2}e^{-2r}+\frac{1}{r}\E(X_{\tau_{2}}-X_{\tau_1})<\infty.
\end{align*}

Now the Borel-Cantelli lemma implies that
$$\P(A_{m} \text{ for infinitely many } m\geq 1)=0.$$
Hence, the walk $(\widetilde S_n)_{n\geq 0}$ almost surely crosses $0$ finitely many times.
Therefore, also
the walk  $(S_n)_{n\geq 0}$ crosses the vertical line $\{0\}\times\R$ finitely many times
and visits just finitely many points on one side of that line.
\qed

\section{Two parallel lines with thinned Poisson processes}
\label{par_indep}

Let $E = \R\times \{0,r\}$ and let $d$ be the
Euclidean distance. 
Moreover, let $0<p\leq 1$.
We define a point process $\Pi$ on $E$ as follows.

Let $\widehat \Pi$ be a homogeneous Poisson process with rate 1.
Let $\Pi^0$ and $\Pi^r$ be two (dependent, in general) thinnings of $\widehat \Pi$ generated as follows.
For all $X\in\widehat \Pi$, do one of the following: 
\begin{itemize}
 \item With probability $1-p$, duplicate the point $X$ and make it a point of both processes $\Pi^0$ and $\Pi^r$.
 \item With probability $p$, assign $X$ to either $\Pi^0$
or $\Pi^r$ (but not both), with probability $1/2$ each.
\end{itemize}

Let now $\Pi$ be the point process on $E$ with 
$\Pi\cap (B\times \{0\})=(\Pi^0\cap B)\times\{0\}$ and $\Pi\cap (B \times \{r\})=(\Pi^r\cap B)\times\{r\}$ for all Borel sets $B\subset \R$.

We study the greedy walk on $\Pi$ defined by \eqref{defwalk} that starts at $S_0=(0,0)$. 
Note that $(0,0)\notin\Pi$ with probability 1.
Sometimes in the proofs we consider the greedy walk starting from another point $x\in E$. 
We emphasise this by writing the superscript $x$ in the sequence $(S^{x}_n)_{n\geq 0}$.

If $p=0$, $\Pi^0$ and $\Pi^r$ are identical. This was studied in Section \ref{parallel_same_process}
where we showed that the greedy walk does not visit all points of the process.
In this section we consider $p>0$ and we get the opposite result:
\begin{theorem}
\label{thm:thinning}
For any $p>0$, the greedy walk visits all points of $\Pi$ almost surely.
\end{theorem}
For $p=1$, the processes $\Pi^0$ and $\Pi^r$ defined above are two independent Poisson processes with rate $1/2$.
When $0<p<1$, the process $\Pi$ has some ``doublë́'' points and  $\Pi^0$ and $\Pi^r$ are not independent.
However, for $p=1$ and $p\in(0,1)$ the behaviour of the walk is similar and, thus, we study these cases together.

Let us now introduce some notation and definitions that we use throughout this section.
We call the projections of elements of $E$ on their first coordinate {\em shadows} of the elements of $E$ on $\R$ and we denote it by $\thinspace\widehat \cdot\thinspace$.
For example, the shadow $\widehat x$ of a point $x=(x_1,x_2)\in E$  is $x_1$.
The point process $\Pi$ is defined from the process $\widehat \Pi$ that can be seen as a shadow of $\Pi$,
that is, $x\in\widehat\Pi$ if and only if there exists $\iota\in\{0,r\}$ such that $(x,\iota)\in \Pi$.
Also, $(\widehat S_n)_{n\geq 0}$ is the shadow of the greedy walk $(S_n)_{n\geq 0}$, that is, $(\widehat S_n)_{n\geq 0}$ contains just the information about the first coordinates of the locations of the walk.

Let $\Pi_n=\Pi\setminus \{S_0,S_1,S_2,\dots,S_{n-1}\}$ be the set of all points of $\Pi$ that are not visited before time $n$.
Let $\Pi^0_n$ and $\Pi^r_n$ be the shadows of $\Pi_n\cap (\R\times\{0\})$ and $\Pi_n\cap(\R\times\{r\})$, respectively. 
Moreover, $\widehat \Pi_n$ denotes the shadow of $\Pi_n$, that is, the set of the first coordinates of the points of $\Pi_n$.

Define the {\em shift operator} $\theta_{(x,\iota)}$, $(x,\iota)\in E$,  on $\Pi$,
by $\theta_{(x,0)}\Pi=\left(\left(\Pi^0-x\right)\times\{0\}\right)\cup \left(\left(\Pi^r-x\right)\times\{r\}\right)$ and 
$\theta_{(x,r)}\Pi=\left(\left(\Pi^r-x\right)\times\{0\}\right)\cup \left(\left(\Pi^0-x\right)\times\{r\}\right)$.
Define also the {\em mirroring operator} $\sigma$ by 
$\sigma \Pi=\left(\left( -\Pi^0\right)\times\{0\}\right)\cup\left(\left( -\Pi^r\right)\times\{r\}\right)$.

For any subset $A$ of $\R$ define
$$T_A=\inf\{n\geq 1: \widehat S_n\in A\}$$
to be the first time the shadow of the greedy walk enters $A$
and write $T_x$ for $T_{\{x\}}$.
Let $T^R_x=\inf\{n\geq 0:\widehat \Pi_{n+1}(x)=0\}$, that is, $T^R_x$ is the time
when both points $(x,0)$ and $(x,r)$ are visited. 
If exactly one of $(x,0)$ and $(x,r)$ is in $\Pi$ then $T^R_x=T_x$.



Note that for any $x>0$ we have $T_{[x,\infty)}<\infty$ or $T_{(-\infty,0)}<\infty$
because $\widehat\Pi[0,x]<\infty$, almost surely, and the walk exits $[0,x]\times\{0,r\}$ in a finite time.
Define now the variable $D_{x}$, $x\in \widehat \Pi$, $x>0$,
as follows.
If $T_{(-\infty,0)}<T_{[x,\infty)}$ then $D_{x}=0$. 
Otherwise, $T_{[x,\infty)}<T_{(-\infty,0)}$ and 
$0< \widehat S_1,\widehat S_2,\dots,\widehat S_{T_{[x,\infty)}-1}<x,\ \widehat S_{T_{[x,\infty)}}\geq x$. 
We can label the remaining points of $\widehat \Pi_{T_{[x,\infty)}}$ in the interval $(0,x)$ by $z_1,z_2,\dots,z_{n-1}$ so that 
\begin{eqnarray}
\label{rem_marks}
 0=z_n<z_{n-1}<\dots<z_1<z_0=x.
\end{eqnarray}
Let then
\begin{eqnarray}
\label{D_x}
 D_x=\max_{0\leq i\leq n-1}\left\{(z_i-z_{i+1})-(x-z_i)\right\}=\max_{0\leq i\leq n-1}\left\{2z_i-z_{i+1}-x\right\}.
\end{eqnarray}
Since $2z_0-z_{1}-x=x-z_1>0$ and $2z_i-z_{i+1}-x\leq 2z_i-x\leq x$ for $0\leq i\leq n-1$, we have
$0<D_x\leq x$.
The variable $D_{x}$  measures how big should be the distance between $x$ and 
the closest point to $x$ in $\widehat \Pi\cap (x,\infty)$,
so that the walk possibly visits a point in $(-\infty, 0)\times\{0,r\}$ 
before visiting any point in $(x,\infty)\times\{0,r\}$.

We prove Theorem \ref{thm:thinning} in a similar way as Rolla {\em et al.} \cite{RST14}
prove that the greedy walk visits all marks attached to the points of a homogeneous Poisson process on $\R$, 
where each point has 
one mark with probability $p$ or two marks with probability $1-p$.
The idea of the proof is the following. 
We define first a subset $\Xi$ of $\Pi^0$, which is stationary and ergodic.
Then, using the definition and properties of $\Xi$, we are able to show
that  there exists $d_0>0$ such that, almost surely, $D_{X}<d_0$ for infinitely many $X\in\widehat \Pi$.   
This we use to show that events
$A_k$, which we define later in \eqref{A_k_events},
occur for infinitely many $k>0$, almost surely.
Then we show that whenever $A_k$ occurs, the greedy walk visits
$(-\infty, 0)\times\{0,r\}$ in a finite time. 
Therefore, we can conclude that $T_{(-\infty,0)}<\infty$, almost surely.
Finally, using repeatedly the fact that $T_{(-\infty,0)}<\infty$, almost surely,
we are able to show that the greedy walk on $\Pi$ crosses the vertical line $\{0\}\times\R$ infinitely many times
and, thus, almost surely visits all points of $\Pi$.

Let us first discuss general properties of the greedy walk on $E$ that do 
not depend on the definition of the point process $\Pi$.
If the greedy walk visits two points $(x,\iota)$ and $(y,\iota)$ on a line $\iota\in\{0,r\}$, 
without changing lines in between those visits, 
then it visits also all the points between $x$ and $y$ on line $\iota$.
So the walk clears some intervals on the lines, but because of changing the lines it  
possibly leaves some unvisited points between those intervals.
The following lemma shows that if the walk omits two points on different lines
between visited intervals, then the horizontal distance between those points is greater than $r$.
Moreover,
looking from a point that is to the right of those two points, 
the closer point is always the point that is more to the right. 
Thus, when the walk returns to a partly visited set, 
it always visits first the rightmost remaining point in this set.

\begin{lemma}
\label{distance}
 Let $a=\min_{0\leq i\leq n} \widehat S_i$ and $b=\max_{0\leq i\leq n} \widehat S_i$.
\begin{itemize}
 \item[(a)] Let 
$a\leq x<y\leq b$
such that $(x,\iota_x),(y,\iota_y)\in \Pi_{n+1}$ and $\iota_x\neq\iota_y$. Then 
$y-x>r$.
\item[(b)] Let $a\leq x<y\leq z\leq b$ such that $(x,\iota_x),(y,\iota_y)\in \Pi_{n+1}$, $(z,\iota_z)\in\Pi_n$.
Then $d((y,\iota_y),(z,\iota_z))<d((x,\iota_x),(z,\iota_z))$.
\end{itemize}
\end{lemma}
\proof 
(a) Let $x,y$ be as in the lemma and suppose on the contrary that $y-x\leq r$.
Let $I_1=\left((a,x)\times\{\iota_x\}\right)\cup\left((a,y)\times\{\iota_y\}\right)$ and $I_2=\left((x,b)\times\{\iota_x\}\right)\cup\left((y,b)\times\{\iota_y\}\right)$.
Since the greedy walk has visited points at $a$ and $b$ up to time $n$, 
but not the points $(x,\iota_x)$ and $(y,\iota_y)$,
the walk before time $n$ moved from $I_1$ to $I_2$ or in the opposite direction.
Because of the assumption $y-x\leq r$, a point in $I_1$ is closer to $(x,\iota_x)$ or $(y,\iota_y)$ than
to any point in $I_2$, and conversely. 
Thus, if $y-x\leq r$ it is impossible to visit a point at $a$ and a point at $b$ without visiting 
either $(x,\iota_x)$ or $(y,\iota_y)$, which contradicts the assumptions of the lemma.

 (b) 
If $\iota_y=\iota_z$ then 
$d((y,\iota_y),(z,\iota_z))=z-y<z-x\leq d((x,\iota_x),(z,\iota_z)).$
If $\iota_y\neq\iota_z$ and $\iota_x=\iota_y$ then 
 $d((y,\iota_y),(z,\iota_z))^2=(z-y)^2+r^2<(z-x)^2+r^2\leq d((x,\iota_x),(z,\iota_z))^2$.
If $\iota_y\neq\iota_z$ and $\iota_x\neq \iota_y$, it follows 
from  part (a) of this lemma that $y-x>r$. 
This yields
\begin{align*}
 d((y,\iota_y),(z,\iota_z))^2&=(z-y)^2+r^2<(z-y+r)^2\\
&<(z-y+y-x)^2=d((x,\iota_x),(z,\iota_z))^2.
\end{align*}
\vskip-8mm\hfill$\Box$\vskip3mm

Let $(c,\iota_c)\in\Pi$, $c>0$, be a point far enough from $(0,0)$.
If $(c,\iota_c)$ is not visited when the walk goes for the first time from $(-\infty,c)\times\{0,r\}$ to $(c,\infty)\times\{0,r\}$, then this point is visited before the walk 
returns to $(-\infty,c)\times\{0,r\}$. 
Furthermore, as we show in the following lemma, when the walk is finally at the location $(c,\iota_c)$, then there is no unvisited points left in $(c,\max_{0\leq i \leq T^R_c}\widehat S_i]\times\{0,r\}$.
If $(c,\iota_c)$ is close to $(0,0)$ and $\iota_c=r$, then this is not always true because the greedy walk might jump several times over $(0,0)$ (and $(c,0)$) without visiting any point on line $r$.

\begin{lemma}
\label{empty_int}
Let $a=\min_{0\leq i\leq n} \widehat S_i$ and $b=\max_{0\leq i\leq n} \widehat S_i$.
 Let $a\leq c\leq b$, $\iota_c\in\{0,r\}$,  such that $(c,\iota_c)\in\Pi_{n+1}$.
If $\widehat S_n\geq c$ and $(c,\iota_c)$
is visited before any other point of $\Pi_n$ in $(-\infty,c)\times\{0,r\}$,
then all points of $\Pi$ in $(c,\max_{0\leq i \leq T^R_c}\widehat S_i]\times\{0,r\}$ are visited before time $T_c^R$.
\end{lemma}
\proof
 Let us denote $\max_{0\leq i \leq T^R_c}\widehat S_i$ by $M$.
It is easy to see that the lemma is true for $c=M$.
Thus, assume that $M>c$ and
let $j$ be the first time 
for which $\widehat S_j=M$. 
Then, for $j+1\leq k\leq T^R_c-1$,  $S_k$ is in $(c,M]\times\{0,r\}$.

Let $x$ be the rightmost point of 
$\widehat \Pi_{j+1}$ in the interval $[c,M]$ and let $\iota_x$ be the corresponding line of that point.
Note that at time $j+1$ there is only one point left at $x$.
If $x=c$, the claim of the lemma follows directly.
Otherwise,
by Lemma \ref{distance} (b),
the point $(x,\iota_x)$ is closer to $S_j$ than any other
point in $[c,x)\times\{0,r\}$. 
Thus, $S_{j+1}=(x,\iota_x)$ and there are no points left in $[x,M]\times\{0,r\}$ after time $j+1$.
Repeating the same arguments we can see that
the walk successively visits all the remaining points 
in $[c,M]\times\{0,r\}$ until it reaches $c$.
Hence, at time $T^R_c$ the set $(c,M]\times\{0,r\}$ is empty.
\qed

To show that $D_X<\infty$ for infinitely many $X\in\widehat \Pi$, the crucial will be the subset $\Xi$ of $\Pi^0$
which we define as follows.
For $X\in\Pi^0$, let $L^0(X)=\max\{Y\in \Pi^0: Y<X\}$. 
Then let
\begin{align*}
\Xi=\{X\in \Pi^0:&~ \ S_n^{(L^0(X),0)}\in \left([L^0(X),\infty)\times\{0,r\}\right)\setminus \{(X,0)\} \text{ for all }n\geq 1 \text{ and } \\
& \ \ S_n^{(L^0(X),r)}\in \left([L^0(X),\infty)\times\{0,r\}\right)\setminus \{(X,0)\} \text{ for all }n\geq 1 \},
\end{align*}
that is, $\Xi$ contains all points $X\in\Pi^0$ such that
if we start a walk from $(L^0(X),0)$ or $(L^0(X),r)$, then the walk stays always in 
$[L^0(X),\infty)\times\{0,r\}$ and it never visits $(X,0)$.
We have defined $\Xi$ in this way because then for every $X\in\Xi$, whenever the greedy walk approaches  
$(X,0)$ and $(X,r)$ from their left, the walk passes around $(X,0)$ and it never comes back to visit $(X,0)$.
Thus, if $\Xi$ is non-empty, there are some points of $\Pi$ that are never visited.
 
The point process $\Pi$ is generated from the homogeneous Poisson process $\widehat \Pi$ and thus it is  stationary and ergodic.
The process $\Xi$ is defined as a function of the points in $\Pi$ and therefore
  $\Xi$ is a stationary and ergodic process in $\R$. 
Therefore, $\Xi$ is almost surely the empty set or it is almost surely a non-empty set, in which case $\Xi$ has a positive rate. 
We look at these two cases in the next two lemmas.

For the first lemma we need the random variable $W_{(x,\iota)}$, $x\in\R$, $\iota\in\{0,r\}$, which measures a sufficient horizontal distance from ${(x,\iota)}$ to the rightmost point in the set
$(-\infty,x)\times\{0,r\}$, so that the greedy walk starting from $(x,\iota)$ never visits that set.
For $x\in\R$ let $\Pi^x=\Pi\cap\left([x,\infty)\times\{0,r\}\right)$.
Consider for the moment the greedy walk on $\Pi^x$ starting from $(x,\iota)$ defined by \eqref{defwalk} and
let 
$$W_{(x,\iota)}=\inf_{i\geq 0}\{\widehat S_i-d(S_i,S_{i+1})\}.$$
Let $L^0(x)=\max\{Y\in\Pi^0:Y<x\}$ and $L^r(x)=\max\{Y\in\Pi^r:Y<x\}$.
If
for some $x\in\R$, $\iota\in\{0,r\}$, $\vert W_{(x,\iota)}\vert<\infty$ and $\max\{L^0(x),L^r(x)\}<W_{(x,\iota)}$,   then for all $i\geq 0$
\begin{align*}
 \min\{d(S_i,(L^0(x),0)),d(S_i,(L^r(x),r))\}&\geq \widehat S_i-\max\{L^0(x),L^r(x)\}\\
&> \widehat S_i-W_{(x,\iota)}\geq d(S_i,S_{i+1}).
\end{align*}
Therefore, the walk on $\Pi$ starting from $(x,\iota)$ coincides with the walk on $\Pi^x$, because for every $i\geq 0$ the point
$S_{i}$ is closer to $S_{i+1}$ than to any point in $\Pi\setminus\Pi^x$.
Thus, the walk does not visit $\Pi\setminus\Pi^x$.
The opposite is also true, i.e.\
 if the walk on $\Pi$ starting from $(x,\iota)$ coincides with the walk on $\Pi^x$,
then $\vert W_{(x,\iota)}\vert<\infty$.

\begin{lemma}
 \label{xi_empty}
If $\Xi$ is almost surely the empty set, then $\Psi=\{(X,\iota)\in\Pi:\vert W_{(X,\iota)}\vert <\infty\}$
is also almost surely the empty set.
\end{lemma}
\proof
For $x\in\R$, let $L^0(x)$ and $L^r(x)$ be as above. 
Since $W_{(X,\iota)}$ is identically distributed for all $(X,\iota)\in\Pi$, 
$\Psi$ is stationary and ergodic.
Suppose, on the contrary, that $\Psi$ is a non-empty set.
For $d>0$ let $\Psi_d=\{(X,r)\in \Psi:\vert W_{(X,r)}\vert<d\}$ and note that $\bigcup_{d> 0}\Psi_d=\Psi\cap (\R\times\{r\})$.
Thus, there exists $d$ large enough so that $\Psi_d$ is a non-empty set which is stationary and ergodic.
Let $\widetilde \Psi_d$ be the set of all $(X,r)\in \Psi_d$ such that the points $(L^0(X),0)$, $(L^r(X),r)$ and $(L^0(L^0(X)),0)$ satisfy the following.
First, these points are in $(-\infty,X-d)\times\{0,r\}$.
Second, the point $(L^r(X),r)$ is closer to $(X,r)$ than to $(L^0(X),0)$.
Third, the greedy walks starting from $(L^0(L^0(X)),0)$ and $(L^0(L^0(X)),r)$
visit only the points in $[L^0(L^0(X)),\infty)\times\{0,r\}$ and these walks visit $(L^r(X),r)$ 
before visiting $(L^0(X),0)$.
Since these three conditions have a positive probability, which is independent of  $W_{(X,r)}$, 
$\widetilde  \Psi_d$ is also almost surely a 
non-empty set.
But, by definition of $\Xi$, for every $(X,r)\in\widetilde \Psi_d$, $L^0(X)\in\Xi$ and $\Xi$ is a non-empty set, which is a contradiction.
\qed

In the following lemma we use the random variable $\widetilde D_X$, $X\in\Pi_0$, which can be compared to $D_X$ defined in \eqref{D_x}. 
For $X\in\Pi^0\setminus\Xi$ set $\widetilde D_X=0$.
For $X\in \Xi$ denote the points of $\Xi\cap(\infty,X]$ in decreasing order
$$\dots<\widetilde z_2<\widetilde z_1<\widetilde z_0=X$$
and define $\widetilde D_X$ as
\begin{align}
\label{tildeDy}
\widetilde D_X =&\sup_{i\geq 0}\left\{2\widetilde z_i-\widetilde z_{i+1}-X\right\}.
\end{align}

For $d_1, d_2>0$ define
\begin{align*}
 \Xi_{d_1,d_2}=\{X\in\Xi:&~
\widetilde D_{X}<d_1 \text{ and there exists }Y\in\widehat \Pi\text{ such that }0<Y-X<d_2\text{ and }\\
&~  \widehat\Pi\cap(Y, Y+r)=\emptyset\},
\end{align*}
i.e.\ $X\in\Xi$ belongs to $\Xi_{d_1,d_2}$  if 
$\widetilde D_{X}<d_1$ and at the distance less than $d_2$ from $X$
there is an interval of length $r$ where there is no points of $\widehat \Pi$.
We consider the empty intervals of length $r$ in $\Pi$, because whenever 
$\widehat \Pi\cap(Y, Y+r)=\emptyset$ for some $Y\in\widehat \Pi$, $Y>0$,
the greedy walk is forced to visit the points  $(Y,0)$ and $(Y,r)$, before crossing the interval and visiting a point in $[Y+r,\infty)\times\{0,r\}$.

\begin{lemma}
\label{d_1_new}
 If $\Xi$ is almost surely a non-empty set, then 
there exist $d_1, d_2>0$ such that $\Xi_{d_1,d_2}$ is almost surely a non-empty set.
Moreover,  $\Xi_{d_1,d_2}$ is a stationary and ergodic process.
\end{lemma}
\proof

If $\Xi$ is non-empty, the rate $\delta$ of $\Xi$ is positive. 
Then for $X\in\Xi$, 
$\lim_{i\rightarrow\infty}  \frac{X-\widetilde z_i}{i} =\delta^{-1}$, almost surely.
Also, $\lim_{i\rightarrow\infty} \frac{\widetilde z_i-\widetilde z_{i+1}}{i}=0$, almost surely.
Hence, for all large $i$, $2\widetilde z_i-\widetilde z_{i+1}-X=
(\widetilde z_i-\widetilde z_{i+1})-(X-\widetilde z_i)<0$ 
and we can conclude that $\widetilde D_X$ is almost surely finite.
Therefore, there exists $d_1$ such that $\{X\in \Pi^0: X\in\Xi,~ \widetilde D_X<d_1\}$ is 
with positive probability a non-empty set.
Moreover, since $\widetilde D_X$ is identically distributed for all $X\in \Xi$,  
$\{X\in\Pi^0:X\in\Xi,~ \widetilde D_X<d_1\}$ 
is a stationary and ergodic process with positive rate.

Almost surely, the gap between two neighbouring points of the homogeneous Poisson process $\widehat\Pi$ 
is  infinitely often greater than $r$. 
Thus, also
$\Pi\cap\left((Y, Y+r)\times \{0,r\}\right)=\emptyset$ for infinitely many $Y\in\widehat\Pi$ and
all such $Y$ form a stationary and ergodic process.
Thus, $\Xi_{d_1,d_2}$ is also stationary and ergodic for all $d_2>0$.
%
Since $\bigcup_{d_2> 0} \Xi_{d_1, d_2}=\{X\in \Pi^0: X\in\Xi, ~\widetilde D_X<d_1\}$ and this is not empty when $d_1$ is large enough, 
we can choose $d_2$ large enough so that $\Xi_{d_1,d_2}$ is almost surely a non-empty set.
\qed

We study the greedy walk starting from the point $(0,0)$, which is, almost surely, not a point of $\Pi$.
From now on denote the points of $\widehat \Pi$ by 
 $$\ldots<X_{-2} <X_{-1}\leq 0< X_1<X_2< \dots $$
and denote the points of $\Pi^0$ and $\Pi^r$ by
 $$\ldots<X_{-2}^0 <X_{-1}^0\leq 0< X_1^0<X_2^0< \dots  \text{and}\ldots<X_{-2}^r <X_{-1}^r\leq 0< X_1^r<X_2^r< \dots,$$
 respectively.

Now we are ready to prove that $D_{X}<d_0$ for infinitely many $X\in\widehat \Pi$. 
We use here the definition of $\Xi$ and divide the proof in two parts, depending weather  $\Xi$ is almost surely the empty set or a non-empty set.

\begin{lemma}
\label{d_0}
 There exists $d_0<\infty$ such that, almost surely, $D_{X_k}<d_0$ and $X_{k+1}-X_k>r$ for infinitely many $k>0$.
\end{lemma}
\proof

Assume first that $\Xi$ is almost surely the empty set.
Then, by Lemma \ref{xi_empty}, 
$\Psi=\{(X,\iota)\in\Pi: \vert W_{(X,\iota)}\vert <\infty\}$
is almost surely the empty set.
Moreover, $\{(X,0)\in\Pi:X\notin \Pi^r, \vert W_{(X,0)}\vert <\infty \}$ is also the empty set.
The greedy walk $(S^{(0,0)}_n)_{n\geq 0}$ on $\Pi\cap \left([0,\infty)\times\{0,r\}\right)$ has the same law as the greedy walk $(S^{(X,0)}_n)_{n\geq 0}$ on $\Pi\cap \left([X,\infty)\times\{0,r\}\right)$ shifted for $(X,0)$, where $X\in\Pi^0, X\notin \Pi^r$.
This implies that $\P(\vert W_{(0,0)}\vert <\infty)=0$
and  the greedy walk $(S^{(0,0)}_n)_{n\geq 0}$ almost surely visits a point in $(-\infty,0)\times\{0,r\}$ in a finite time. 
Then it follows from the definition of $D_x$ that $D_x=0$ for all large enough $x> 0$. 
Since $X_{k+1}-X_k>r$ for infinitely many $k>0$, the 
claim of the lemma holds for any $d_0>0$.

Assume now that $\Xi$ is not empty. 
Then, by Lemma \ref{d_1_new}, we can find $d_1$ and $d_2$ large enough so that $\Xi_{d_1,d_2}$ is a non-empty set.
We first show that $D_X\leq d_1$ for infinitely many $X\in \Xi$, $X>0$, and then we prove that
$D_{X_k}<d_1+d_2$ and $X_{k+1}-X_k>r$ for infinitely many $k>0$.

For $X\in\Xi$, $X>0$, denote the points
of $\Xi$ in $[0,X]$ by $\widetilde z_0,\widetilde z_1,\dots,\widetilde z_{\widetilde n-1}$ so that 
$$0=\widetilde z_{\widetilde n}<\widetilde z_{\widetilde n-1}<\dots<\widetilde z_1<\widetilde z_0=X$$
and define $\widetilde D_X'$, an analogue of $D_X$ and a restricted version of $\widetilde D_X$, as
\begin{align}
\label{tildeDy2}
\widetilde D_X'=\max_{0\leq  i \leq {\widetilde n-1}}\{2\widetilde z_i-\widetilde z_{i+1}-X\}=
\max\left\{\max_{0\leq  i \leq {\widetilde n-2}}\{2\widetilde z_i-\widetilde z_{i+1}-X\},2\widetilde z_{\widetilde n-1}-X\right\}.
\end{align}

Let $\xi=\min\{Y\in\Xi:~Y>0\}$ and note that $z_{\widetilde n-1}=\xi$ for every $X\in \Xi$, $X>0$.
Then for $X\in \Xi$, $X>2\xi$ we have $2z_{\widetilde n-1}-X=2\xi-X<0$.
Since $\widetilde D_X'\geq 2\widetilde z_0-\widetilde z_{1}-X=X-\widetilde z_{1}>0$,
the term $2z_{\widetilde n-1}-X$ does not contribute to $\widetilde D_X'$.
From the definition of $\Xi_{d_1,d_2}$, we have $\widetilde D_X<d_1$ for infinitely many $X\in\Xi$, almost surely.
When $X>2\xi$, $\widetilde D_X'$ is the maximum of a finite subset of the values in \eqref{tildeDy}
and thus $\widetilde D_X'\leq \widetilde D_X<d_1$ for infinitely many $X\in\Xi$, $X>2\xi$.

We prove now that $D_X\leq \widetilde D_X'$ in two steps. 
First, we show that the points used in the definition of $\widetilde D_X'$ are a subset of the points used in the definition of $D_X$.
Second, we show that adding a new point to the definition \eqref{tildeDy2} decreases the value of the maximum.
Let $\xi\in\Xi$, $\xi>0$. 
Before visiting any point in $[\xi,\infty)\times\{0,r\}$, the greedy walk starting from $(0,0)$ visits 
the leftmost point 
on one of the lines in $[L^0(\xi),\infty)\times\{0,r\}$, where $L^0(\xi)=\max\{Y\in\Pi^0:Y<\xi\}$.
That is, the greedy walk visits $(L^0(\xi),0)$ or the closest point to the right of $(L^0(\xi),r)$ on the line $r$ 
(which is $(L^0(\xi),r)$ if that point exists).
By the definition of $\Xi$, the greedy walk starting from one of these two points never visits
$(\xi,0)$ and it never visits any point to the left of $\{L^0(\xi)\}\times\{0,r\}$.
Therefore, 
once the greedy walk starting from $(0,0)$ enters $[L^0(\xi),\infty)\times\{0,r\}$, it continues on the path of
one of these two walks.
Thus, the greedy walk starting from $(0,0)$ does not visit $(\xi,0)$.
From this we can conclude that all points in $\{\Xi\cap(0,\infty)\}\times\{0\}$ are not visited by the greedy walk and 
for $X\in \Xi\cap (0,\infty)$ we have $\{\widetilde z_0,\widetilde z_1,\dots,\widetilde z_{\widetilde n}\}\subset \{z_0,z_1,\dots,z_n\}$,
where $z_0,z_1,\dots,z_n$ are as in \eqref{rem_marks}.

Let  $Y\in \{z_0,z_1,\dots,z_n\}\setminus \{\widetilde z_0,\widetilde z_1,\dots,\widetilde z_{\widetilde n}\}$ and 
find  $j$ such that $\widetilde z_j<Y<\widetilde z_{j-1}$. 
Adding $Y$ to the set $\{\widetilde z_0,\widetilde z_1,\dots,\widetilde z_{\widetilde n}\}$ in the definition \eqref{tildeDy2},
removes
the value $2\widetilde z_{j-1}-\widetilde z_j-X$ 
and adds the values
  $2\widetilde z_{j-1}-Y-X$ and $2Y-\widetilde z_j-X$.
Since, $2\widetilde z_{j-1}-Y-X<2\widetilde z_{j-1}-\widetilde z_j-X$ and 
$2Y-\widetilde z_j-X<2\widetilde z_{j-1}-\widetilde z_j-X$, 
the point at $Y$ added to $\{\widetilde z_0,\widetilde z_1,\dots,\widetilde z_{\widetilde n}\}$ 
decreases the value of $\widetilde D_X'$ or leaves it unchanged.
Since both $D_X$ and $\widetilde D_X'$ are defined on $\{\widetilde z_0,\widetilde z_1,\dots,\widetilde z_{\widetilde n}\}$,
but $D_X$ has also points $\{z_0,z_1,\dots,z_n\}\setminus \{\widetilde z_0,\widetilde z_1,\dots,\widetilde z_{\widetilde n}\}$
in the definition,
we can conclude that
$D_X\leq \widetilde D_X'$ for all $X\in \Xi$, $X>0$.
Thus $D_X<d_1$ for infinitely many $X\in \Xi$, $X>0$.

This together with Lemma \ref{d_1_new} implies that for infinitely many 
$X\in\Xi \cap (0,\infty)$, $D_X<d_1$ and there exists $Y\in \widehat \Pi$ such that $0< Y-X<d_2$
and $\Pi(( Y, Y+r)\times \{0,r\})=0$.
Choose one such $X$ and let $k$ be such that $X_k-X<d_2$ and $X_{k+1}-X_k>r$.
If $T_{(-\infty,0)}<T_{[X_k,\infty)}$, then $D_{X_k}=0<d_1+d_2$.
Otherwise, $T_{[X_k,+\infty)}<T_{(-\infty,0)}$,  and we can denote the points of $\widehat \Pi_{T_{[X_k,\infty)}}$ 
in  $(0,X_k)$ by $z_1,z_2,\dots,z_{n-1}$
so that $0=z_n<z_{n-1}<\dots<z_1<z_0=X_k$.
By the definition of $\Xi$, the point $(X,0)$ is never visited by the walk and 
thus there exists $j$ such that $z_j=X$.
Then we have
\begin{align*}
 D_{X_k}&=\max_{0\leq i\leq n-1}\left\{2 z_i- z_{i+1}-X_k\right\}\\
&\leq \max\{ \max_{j\leq i\leq n-1}\left\{2 z_i- z_{i+1}-X\right\}-(X_k-X),\max_{0\leq i\leq j-1}\left\{2 z_i- z_{i+1}-X_k\right\}\}\\
&\leq \max\{D_{X},X_k-X\}\leq  D_{X}+X_k-X< d_1+d_2,
\end{align*}
where in the second inequality we use the fact that, by the definition of the points in $\Xi$,
the walk does not visit any point in $(0,X)\times\{0,r\}$ 
after time $T_{[X,\infty)}$
and therefore $D_{X}=\max_{j\leq i\leq n-1}\left\{2 z_i- z_{i+1}-X\right\}$.

 Let now $d_0=d_1+d_2$.
 Since there are, almost surely, infinitely many $X\in \Xi$ and $k$
 such that $D_{X}<d_1$, $X_k-X<d_2$ and $X_{k+1}-X_k>r$,
 it follows that  $D_{X_k}<d_0$ and $X_{k+1}-X_k>r$ for infinitely many $k>0$,  almost surely, which proves the claim of the lemma.
\qed

Since $D_{X_k}<d_0$ for infinitely many $k$, almost surely, 
one should expect that also $X_{k+1}-X_k>d_0>D_{X_k}$ for infinitely many $k$. 
That is exactly what we show next, but let us first state the extended Borel--Cantelli Lemma
which we use in the proof.

\begin{lemma}[Extended Borel--Cantelli lemma, {\cite[Corollary 6.20]{Ka}}]
\label{extended}
Let $\mathcal{F}_n$, $n\geq 0$, be a filtration with $\mathcal{F}_0=\{0,\Omega\}$
and let $A_n\in\mathcal{F}_n,n\geq 1$.  Then a.s.
\begin{align*}
 \{A_n \text{ i.o.}\}=\left\{\sum_{n=1}^{\infty}\mathbb{P}[A_n~\vert~\mathcal{F}_{n-1}]=\infty\right\}.
\end{align*}
\end{lemma}

\begin{lemma}
\label{prob_is_1}
Almost surely, the events 
\begin{align}
\label{A_k_events}
 A_k=\{X_{k+1}-X_k>D_{X_k}-X_{-1}+r\}
\end{align}
occur for infinitely many $k>0$.
\end{lemma}
\proof
For $k>0$ let $j_k^0=\max\{i:X_i^0\leq X_k\}$ and $j_k^r=\max\{i:X_i^r\leq X_k\}$.
Furthermore, define the $\sigma$-algebra 
$$\mathcal F_k=\sigma((X_{-1}^0,0),(X_0^0,0),\dots,(X_{j_k^0}^0,0),(X_{-1}^r,r),(X_{0}^r,r),(X_1^r,r),\dots,(X_{j_k^r}^r,r))$$
and denote by $T_A^{\sigma}$ and $D_x^{\sigma}$ analogues of $T_A$ and $D_x$ for 
the greedy walk on the set of points which generates $\mathcal F_k$.
Assume $T_{[X_k,\infty)}<T_{(-\infty,0)}$.
Then, the greedy walk on $\Pi$ and the walk on the restricted set are the same until time $T_{[X_k,\infty)}$.
If $X_{k+1}-X_k>r$ then $\widehat S_{T_{[X_k,\infty)}}=X_k$, $T_{X_k}^{\sigma}= T_{[X_k,\infty)}$ and
 $D_{X_k}^{\sigma}= D_{X_k}$.

Let $A_k^{\sigma}=\{X_{k+1}-X_k>D_{X_k}^{\sigma}-X_{-1}+r\}$
and observe that $A_k^{\sigma}\in \mathcal F_{k+1}$.
For $d_0>0$ we have 
\begin{eqnarray*}
 \P\left(A_k^{\sigma} ~\vert~\mathcal F_k\right)&\geq& \P\left(D_{X_k}^{\sigma}<d_0,X_{k+1}-X_k>d_0-X_{-1}+r ~\vert~\mathcal F_k\right)\\
&=&\1_{\{D_{X_k}^{\sigma}<d_0\}}\P\left(X_{k+1}-X_k>d_0-X_{-1}+r ~\vert~\mathcal F_k\right)\\
&=& \1_{\{D_{X_k}^{\sigma}<d_0\}}e^{-(d_0-X_{-1}+r)} \quad \text{a.s.}
\end{eqnarray*}
The first equality above holds because $\{D_{X_k}^{\sigma}<d_0\}\in \mathcal F_k$.
The second equality follows from the facts that $X_{-1}\in \mathcal F_k$ and $X_{k+1}-X_k$ is exponentially distributed with mean 1
and independent of $\mathcal F_k$.
By Lemma \ref{d_0}, there exists $d_0$ such that  $D_{X_k}<d_0$ and $X_{k+1}-X_k>r$ 
for infinitely many $k$, almost surely. Since, $D_{X_k}^{\sigma}= D_{X_k}$ whenever $X_{k+1}-X_k>r$, 
also $D^{\sigma}_{X_k}<d_0$ for infinitely many $k$ and, thus, 
$$\sum_{k=1}^{\infty}\P\left(A_k^{\sigma}~\vert~\mathcal F_k\right)=\infty\quad \text{a.s.}$$ 
It follows now from the extended Borel-Cantelli lemma (Lemma \ref{extended}) that 
$$\P(A_k^{\sigma} \text{ for infinitely many } k\geq 1)=1.$$
Since $A_k^{\sigma}\subset \{X_{k+1}-X_k>r\}$ and $A_k=A_k^{\sigma}$ whenever $X_{k+1}-X_k>r$, also
$$\P(A_k \text{ for infinitely many } k\geq 1)=1.$$
\vskip-8mm\hfill$\Box$\vskip3mm

Whenever $A_k$ occurs, as we show in the next lemma, the greedy walk is forced to
visit $(-\infty,0)\times\{0,r\}$ before visiting $[X_{k+1},\infty)\times\{0,r\}$.
Note that the arguments do not depend on the definition of the point process $\Pi$.

\begin{lemma}
\label{povratak}
Almost surely, $T_{(-\infty,0)}<\infty$.
\end{lemma}
\proof
By Lemma \ref{prob_is_1}, the events $A_k=\{X_{k+1}-X_k>D_{X_k}-X_{-1}+r\}$ occur for infinitely many $k$, 
almost surely.
To prove the lemma,
it suffices to prove that whenever $A_k$ occurs, 
then $T_{(-\infty,0)}<T_{[X_{k+1},\infty)}$. 
Because the walk exits $[0,X_{k+1}]\times\{0,r\}$ in, almost surely, finite time,
it follows that $T_{(-\infty,0)}<\infty$, almost surely.

Assume that $T_{[X_{k},\infty)}<T_{(-\infty,0)}$ and $A_k$ occurs for some $k$.
Then $X_{k+1}-X_k>D_{X_k}-X_{-1}+r>r$ and a point in $(0,X_k)\times\{0,r\}$
is closer to a point at $X_k$ than to any point in 
$[X_{k+1},\infty)\times\{0,r\}$. 
Hence, $\widehat S_{T_{[X_{k},\infty)}}=X_k$.

Denote the remaining points of $\widehat\Pi_{T_{[X_{k},\infty)}}$ in the interval $[0,X_{k}]$ 
as in \eqref{rem_marks}.
Note that at time $T_{[X_{k},\infty)}$ there is exactly one unvisited point left at each position $z_1,z_2,\dots,z_{n-1}$  and 
denote by $\iota_1,\iota_2,\dots,\iota_{n-1}$ the corresponding lines of these points.
If there is only one point at $X_k$, let $\iota_0$ be the line of this point.
If there are two points at $X_k$, let $\iota_0$ be the line of the point that is not visited at time $T_{X_{k}}$. 
The point $S_{T_{X_k}}$ is closer to the second point at $X_k$, if such exists, than to any point in $(X_{k+1},\infty)\times\{0,r\}$ because $X_{k+1}-X_k>r$
or to any of the remaining points with shadows at $z_{1}, z_2,\dots, z_{n-1}$ because of Lemma \ref{distance} (b). 
%

For $i=0,1,\dots,n-2$,  from the definition of $D_{X_k}$ \eqref{D_x} we have $D_{X_k}\geq 2z_i-z_{i+1}-X_k$ and 
\begin{align*}
 d((z_i,\iota_i),(z_{i+1},\iota_{i+1}))^2&\leq (z_i-z_{i+1})^2+r^2
< (z_i-z_{i+1}+r)^2 \leq (D_{X_k}+X_{k}-z_i+r)^2\\
&< (D_{X_k}-X_{-1}+r+X_{k}-z_i)^2< (X_{k+1}-z_i)^2\\
&\leq d((z_i,\iota_i),(X_{k+1},\iota_i))^2.
\end{align*}
Thus the point $(z_i,\iota_i)$ is closer to $(z_{i+1},\iota_{i+1})$
than to any point in $[X_{k+1},\infty)\times\{0,r\}$.
Moreover, by Lemma \ref{distance} (b), $(z_{i+1},\iota_{i+1})$ is closer to $(z_i,\iota_i)$
than any other point in $[0,z_{i+1})\times\{0,r\}$. 
Thus, when the walk is at $(z_i,\iota_i)$ it visits  $(z_{i+1},\iota_{i+1})$ next, 
except if $z_i<r$ and
there is a point at $(-\infty,0)\times \{0\}$ which is closer.
In the latter case we have $T_{(-\infty,0)}<T_{[X_{k+1},\infty)}$.

Assume that the walk visits successively the points at $z_1,z_2,\dots, z_{n-1}$.
Hence, all points in 
$(z_{n-1},X_{k+1})\times\{0,r\}$ are visited.
When the greedy walk is at $(z_{n-1},\iota_{n-1})$, the closest unvisited point is 
in $(-\infty,0)\times\{0,r\}$, 
because a point with shadow $X_{-1}$ is closer to $(z_{n-1},\iota_{n-1})$ than any point in $[X_{k+1},\infty)\times\{0,r\}$.
This follows from  
\begin{align*}
d((z_{n-1},\iota_{n-1}),(X_{-1},r-\iota_{n-1}))^2&=(z_{n-1}-X_{-1})^2+r^2\\
&\leq (D_{X_k}+X_k-z_{n-1}-X_{-1})^2+r^2\\
&< (D_{X_k}+X_k-z_{n-1}-X_{-1}+r)^2< (X_{k+1}-z_{n-1})^2\\
&=d((z_{n-1},\iota_{n-1}),(X_{k+1},\iota_{n-1}))^2,
\end{align*}
where in the first inequality above we use that, by the definition of $D_{X_k}$ \eqref{D_x}, $D_{X_k}\geq 2z_{n-1}-X_k$.
Thus, the walk visits $(-\infty,0)\times\{0,r\}$ next. Therefore, also now $T_{(-\infty,0)}<T_{[X_{k+1},\infty)}$, 
which completes the proof.
\qed

Now we are ready to prove Theorem \ref{thm:thinning} where we
use Lemma \ref{povratak} repeatedly to show that the greedy walk
crosses the vertical line $\{0\}\times\R$ infinitely often. 
Therefore the greedy walk visits all the points of $\Pi$.

\proof[Proof of Theorem \ref{thm:thinning}]
From Lemma \ref{povratak} we have $\P(T_{(-\infty,0)}<\infty)=1$, which is equivalent to
\begin{equation}
\label{pr1_changed}
\P\left(T_{(-\infty,0)}<\infty\ \vert\ X_{-1}^0, X_{-1}^r,X_1^0,X_1^r\right)=1\quad \text{a.s.}
\end{equation}
Moreover,  the conditional probability \eqref{pr1_changed} is almost surely $1$ for any absolutely continuous distribution of
$(X_{-1}^0, X_{-1}^r,X_1^0,X_1^r)$ on $(-\infty,0)^2\times(0,\infty)^2$, 
which is independent of
$\Pi\cap\left(\left((X_1^0,\infty)\times \{0\}\right)\cup \big((X_1^r,\infty) \times\{r\}\big)\right)$.

Let $T_0=0$ and
let $T_i$, $i\geq 1$, be the first time the greedy walk
visits a point in $(-\infty,\min_{0\leq n\leq T_{i-1}}\widehat S_n)\times\{0,r\}$ for $i$ odd and
in $(\max_{0\leq n\leq T_{i-1}}\widehat S_n,\infty)\times\{0,r\}$ for $i$ even.
That is, for $i$ odd (even) $T_i$ is the time when the walk visits 
the part of $\Pi$ on the left (right) of $\{0\}\times\R$
 which is unvisited up to time $T_{i-1}$. 
We prove first that $T_i$ is almost surely finite for all $i$
and, thus, the greedy walk, almost surely, crosses the vertical line $\{0\}\times\R$ infinitely many times.
Then we show that it is not possible that a point of $\Pi$ is never visited, and thus the 
walk almost surely visits all points of $\Pi$.

Assume that $T_i$ is finite for some even $i$. 
Let $Y_0=\widehat S_{T_i}$, 
$Y_1^{0}=\min\{Y\in\Pi^0:Y>Y_0\}$ and $Y_1^{r}=\min\{Y\in\Pi^r:Y>Y_0'\}$.
Furthermore, let
$Y_{-1}^{0}=\max \{ Y\in\Pi^0: Y<\min_{0\leq n\leq T_{i}} \widehat S_n\}$ and
$Y_{-1}^{r}=\max \{ Y\in\Pi^r: Y<\min_{0\leq n\leq T_{i}} \widehat S_n\}$.
By the definition of $Y_{-1}^{0}$ and $Y_{-1}^{r}$, at time $T_i$
the set $I_{1}=(-\infty,Y_{-1}^{0}]\times\{0\}\cup(-\infty,Y_{-1}^{r}]\times\{r\}$ is not yet visited.
Also, by definition $Y_1^{0}, Y_1^{r} >\widehat S_{T_i}$ and
the set $I_2=\left([Y_1^{0},\infty)\times\{0\}\right)\cup \big([Y_1^{r},\infty)\times\{r\}\big)$ is never visited before time $T_i$.
Moreover, 
because of the strong Markov property, the distributions of  
$\Pi\cap I_1$ and $\Pi\cap I_2$ 
are independent of the points of $\Pi$ outside $I_1$ and $I_2$.

Let $\Pi'=\Pi\cap \left( I_1 \cup I_2 \right)$.
From \eqref{pr1_changed} we have that the greedy walk on $ \theta _{S_{T_i}}\Pi'$ starting from $(0,0)$,
visits the set $(-\infty,0)\times \{0,r\}$ in a finite time, almost surely.
In other words, the greedy walk on $\Pi'$ 
starting from $S_{T_i}$, almost surely, visits a point in
$I_1$ at some time $T_{i+1}'<\infty$.
The greedy walk on $\Pi_{T_i}$, starting from $S_{T_i}$,
might differ from the walk on $\Pi'$ 
if there are some points outside $I_1$ and $I_2$ that are not visited up to time $T_i$.
Denote the shadows of these points by $c_1,c_2,\dots,c_j$, so that
$\widehat S_{T_i}\geq c_1>c_2>\dots>c_j> \max\{Y_{-1}^{0}, Y_{-1}^{r}\}$.
Because of 
Lemma \ref{distance} (b), a point in $(\widehat S_{T_i},\infty)\times\{0,r\}$ is closer to the point at $c_1$,
than to any of the points at $c_2,c_3,\dots,c_j$.
Hence the walk visits the point at $c_1$ before visiting any of the points at $c_2,c_3,\dots,c_j$.

Let $T_{i+1}^{1} =\min\{T_{c_1}^R, T_{\widehat I_1}\}$, where $\widehat I_1=(-\infty,\max\{Y_{-1}^{0'},Y_{-1}^{r'}\})$.
The greedy walks on $\Pi_{T_i}$ and $\Pi'$ starting from $S_{T_i}$ are the same until the time $T_{i+1}^{1}$.
If $T_{i+1}^{1}=T_{\widehat I_1}$ then $T_{i+1}'=T_{i+1}^{1}=T_{i+1}$ and, thus, $T_{i+1}$ is 
finite.
Otherwise, the point at $c_1$ is visited before any point in $I_1$, so 
$T_{c_1}^R<T_{i+1}'$ and $T_{c_1}$ is, almost surely, finite.
In this case, similarly as above, let us define
 $Y_0^{c_1}=c_1$, 
$Y_1^{0,c_1}=\min \{ Y\in\Pi^0: Y>\max_{0\leq n\leq T_{c_1}} \widehat S_n\}$
and $Y_1^{r,c_1}=\min \{ Y\in\Pi^r: Y>\max_{0\leq n\leq T_{c_1}} \widehat S_n\}$.
Moreover, let $I_1^{c_1}=I_1$ and let 
$I_2^{c_2}=\left([Y_1^{0,c_1},\infty)\times\{0\}\right)\cup \Big([Y_1^{r,c_1},\infty)\times\{r\}\Big)$. 
From Lemma \ref{empty_int} we can deduce that
$(c_1,\max_{0\leq n\leq T_{c_1}^R}\widehat S_n]\times\{0,r\}$
is empty at time $T_{c_1}^R$ and, therefore, $I_2^{c_1}$ contains all points of $\Pi_{T_{c_1}^R}$ to the right of $\{(c_1,0),(c_1,r)\}$. 
Now, by the same arguments as above, it follows that
the walk starting from the point at $c_1$ visits almost surely 
a point in $I_1$ or a point at $c_2$ in a finite time.
If the walk, starting from a point at $c_k$, $1\leq k\leq j-1$, visits a point at $c_{k+1}$ before $I_1$, 
we repeat the same procedure.
Since there are only finitely many such points, the walk 
in almost surely finite time
eventually visits $I_1$ and thus $T_{i+1}<\infty$,
almost surely.

When $i$ is odd, 
we can look at the walk on $\sigma \Pi_{T_{i}}$ starting from $\sigma S_{T_i}$.
Then the same procedure as above
yields $T_{i+1}<\infty$, almost surely.
Therefore, $T_i$ is almost surely finite for all $i$.
Assume now that the walk does not visit all points of $\Pi$ and let $(x,\iota_x)$ 
be a point of $\Pi$ that is never visited.
Then, there is $i_0$ even, such that $x<\widehat S_{T_{i_0}}$ and  
$x-\min_{0\leq n\leq T_{i_0}} \widehat S_n>r$.
Then for all $n\in (T_{i_0},T_{i_0+1}-1)$ such that $\widehat S_n\geq x$,
by the choice of $i_0$, 
$S_n$ is closer to $(x,\iota_x)$ than to any  point in 
$(-\infty,\min_{0\leq n\leq T_{i_0}} \widehat S_n)\times\{0,r\}$.
Also, by Lemma \ref{distance} (b),
$S_n$ is closer to $(x,\iota_x)$ than to any  remaining point in 
$[\min_{0\leq n\leq T_{i_0}} \widehat S_n,x)$.
Hence, the greedy walk visits $(x,\iota_x)$
before time $T_{i_0+1}$, which is a contradiction. 
\qed

\section{Two parallel lines with shifted Poisson processes} 
\label{shift}

Let $E = \R\times \{0,r\}$ and let $d$ be the
Euclidean distance.
We define a point process $\Pi$ on $E$ in the following way.
Let $\Pi^0$ be a homogeneous Poisson process on $\R$ with rate $1$ 
 and let $\Pi^r$ be a copy of $\Pi^0$ shifted by $s$, $0<\vert s\vert <r/\sqrt{3}$, 
i.e.\ $\Pi^r= \{x:x-s\in\Pi^0\}.$ 
Then, let $\Pi$ be a point process on $E$ with $\Pi\cap (B\times \{0\})=(\Pi^0\cap B)\times \{0\}$ 
and $\Pi\cap(B\times\{r\})=(\Pi^r\cap B)\times \{r\}$ for all Borel sets $B\subset \R$.
We consider the greedy walk on $\Pi$ defined by  \eqref{defwalk} starting from $S_0=(0,0)$. 

We call the pair of points $(Y,0)$ and $(Y+s,r)$ {\em shifted copies}.
Moreover, we say that a point of $\Pi$ is an {\em indented point} if it is further away from the vertical line $\{0\}\times \R$ than its shifted copy.
That is, for $s>0$, the indented points are in $(-\infty,0)\times\{0\}$ and $(0,\infty)\times\{r\}$.
For $s<0$, the indented points are in $(0,\infty)\times\{0\}$ and 
$(-\infty,0)\times\{r\}$.

We can divide the points of $\Pi$ into clusters in the following way.
Any two successive points on line $0$ are in the same cluster if their distance is less 
than $\sqrt{r^2+s^2}$, otherwise they are in different clusters.
Moreover, any two points on line $0$  
are in the same cluster only if all points
between those two points belong to that cluster.
The points on line $r$ belong to the cluster of its shifted copy.
Throughout this section, we will call
the closest point to the vertical line $\{0\}\times\R$ of a cluster on each line the {\em leading point} of the cluster. 
Every cluster has one leading indented and one leading unindented point,
except the cluster around $(0,0)$ that has possibly points 
in both $(-\infty,0)\times\{0,r\}$ and $(0,\infty)\times\{0,r\}$.

\begin{figure}[t]
    \centering
    \includegraphics[width=\textwidth]{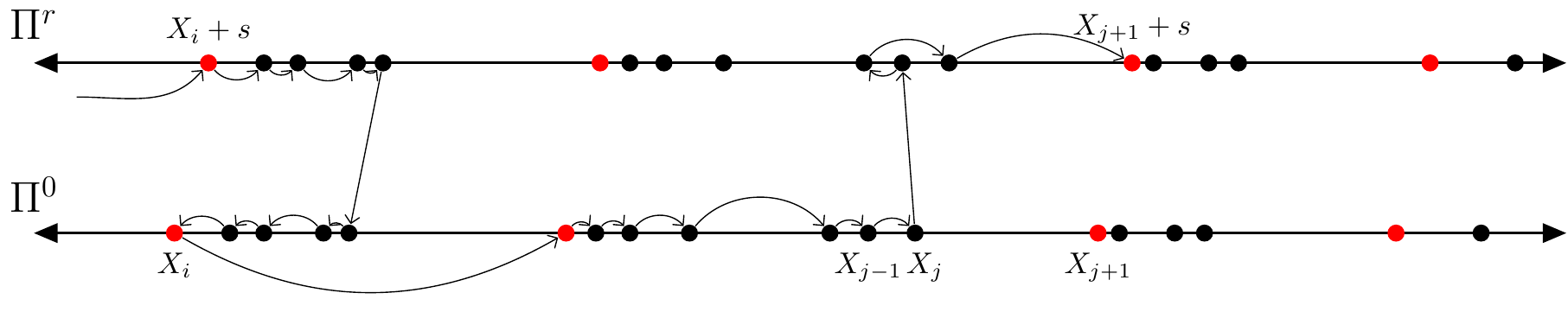}
    \caption{
The leading points of the clusters are marked in gray.
After visiting  the indented leading point $(X_i+s,r)$, 
the walk visits successively all points of its cluster.
This is not always the case when the first visited point of the cluster is the unindented leading point.
In the example above, the point $(X_{j},0)$ is closer to $(X_{j-1}+s,r)$ than to $(X_j+s,r)$, so the walk
moves from  $(X_j,0)$ to $(X_{j-1}+s,r)$. 
Later, when the walk is at $(X_j+s,r)$,  
the closest unvisited point is $(X_{j+1}+s,r)$, which is in a new cluster. 
The greedy walk moves to the next cluster before visiting all points of the current cluster.
}
    \label{shifted_example}
\end{figure}

In Section \ref{parallel_same_process} the points were divided into clusters in a similar way
and we observed that the walk always visits all points of a cluster 
before moving to a new cluster.
This is not the case here.
See Figure \ref{shifted_example} 
for an example where the greedy walk moves to a new cluster before visiting all points in a current cluster.
The points that are not visited during the first visit of a cluster
cause the walk to jump over the vertical line $\{0\}\times\R$ infinitely many times.
Therefore, we obtain here the same result as in Section \ref{par_indep}:
\begin{theorem}
\label{thm_emptyset4}
Almost surely, the greedy walk visits all points of $\Pi$.
\end{theorem}

The proof follows similarly as the proof of Theorem \ref{thm:thinning}.
We change here the definition of the set $\Xi$. 
In addition, the arguments in the first part of the proof of Lemma \ref{d_0_shift}, 
 where we show that if $\Xi$ is almost surely the empty set then the walk almost surely jumps over the starting point,
are different from those in Lemma \ref{d_0}.
Furthermore, in the proof of Theorem \ref{thm_emptyset4} we use the fact that
whenever the walk enters a cluster at its leading unindented point,
then the walk always visits successively  all
points of the cluster. 
This can be explained as follows:
Let $(X,r)$ be the leading indented point and assume that 
$(X,r)$ is the first point that the greedy walk visits in its cluster.
(The case when the walk first visits the leading indented point on the line $0$ can be handled in the same way.)
Denote the closest point to $(X,r)$ on line $r$ by $(Y,r)$ and 
let $a=\vert X-Y\vert$ be the distance between points $(X,r)$ and $(Y,r)$.
The distance between  $(X,r)$ and $(X-s,0)$ is $\sqrt{r^2+s^2}$
and the distance between $(X,r)$ and $(Y-s,0)$ is $\sqrt{r^2+(a-s)^2}$.
Because of the choice $s<r/\sqrt{3}$, we have $2s<\sqrt{r^2+s^2}$. Therefore, if $r\leq a<\sqrt{r^2+s^2}$, 
then 
\begin{align*}
 r^2+(a-s)^2&=a^2+r^2+s^2-2sa>a^2+r^2+s^2-2s\sqrt{r^2+s^2}\\
&=a^2+\sqrt{r^2+s^2}\left(\sqrt{r^2+s^2}-2s\right)>a^2
\end{align*}
Also, if $a<r$, then $r^2+(a-s)^2\geq r^2>a^2$. 
Thus, when $a<\sqrt{r^2+s^2}$ the point $(Y,r)$ is closer to $(X,r)$
than the copies of those two points on line $0$ and thus $(Y,r)$ is visited next.
We can argue in the same way for all the points in this cluster on line $r$, until the walk reaches the outermost point.
The distance from the outermost point of the cluster to the closest point on line $r$ is greater than  $\sqrt{r^2+s^2}$
and the closest unvisited point is its shifted copy.
Once the walk is on line $0$, it visits all remaining points of the cluster,
because the distances between the successive points in the cluster are less than $\sqrt{r^2+s^2}$,
all points of the cluster on line $r$ are already visited and
the distance to any point in another cluster is greater than $\sqrt{r^2+s^2}$.

We define now the set $\Xi$ in a slightly different way than in Section \ref{par_indep}.
For $x\in\R$, let $L^0(x)=\max\{Y\in \Pi^0: Y<x\}$.
Then define 
\begin{align*}
 \Xi=s&+ \left\{X\in \Pi^0:~  S_n^{(X,0)}\in \left((X,\infty)\times\{0,r\}\right)\setminus\{(X+s,r)\} \text{ for all }n\geq 1  \vphantom{\frac{r^2+s^2}{2s}} \right.\\  
&  \left. \text{ and } X-L^0(X)>\frac{r^2+s^2}{2s} \right\},
\end{align*}
that is, $\Xi$ contains all points $X+s\in\Pi^r$ such that the distance between $(X,0)$ and the closest point on the line $0$ to the right of $(X,0)$ is greater than $\frac{r^2+s^2}{2s}$
and if we start a walk from $(X,0)$, then the walk always stays in $[X,\infty)\times\{0,r\}$, but it never visits $(X+s,r)$.
If $X-L^0(X)>\frac{r^2+s^2}{2s}$ and the walk approaches $(X,0)$ and $(X+s,r)$ from their right,
then the walk visits the point $(X,0)$ before visiting $(X+s,r)$. 
Hence, if $X+s\in\Xi$ and $X>0$, then
the point $(X+s,r)$  is never going to be visited.

Note that the set $\Xi$ is a function of a homogeneous Poisson process and hence it is stationary and ergodic.

Let us define now the random variable $W_{x}$, $x\in\R$, that corresponds to the random variable $W_{x,\iota}$ from Section \ref{par_indep}.
For $x\in\R$, let
$$\Pi^x=\Pi\cap\left(\left([x,\infty)\times\{0\}\right)\cup\left([x+s,\infty)\times\{r\}\right)\right).$$
Consider for the moment the greedy walk on $\Pi^x$ starting from $(x,0)$ defined by \eqref{defwalk} and
let 
$$W_x=\inf_{i\geq 0}\{\widehat S_i-d(S_i,S_{i+1})\}.$$
If $s>0$ and if  
for some $x\in\R$ we have $\vert W_x\vert<\infty$ and $L^0(x)+s<W_x$, then for the walk starting from $(x,0)$ and for any $i\geq 0$
it holds 
\begin{align*}
 \min\{d(S_i,(L^0(x),0)),d(S_i,(L^0(x)+s,r))\}\geq \widehat S_i-L^0(x)-s> \widehat S_i-W_x\geq d(S_i,S_{i+1}),
\end{align*}
that is, the point
$S_{i}$ is closer to $S_{i+1}$ than to any point in $\Pi\setminus\Pi^x$.
Therefore, 
the walk on $\Pi$ starting from $(x,0)$ coincides with the walk on
 $\Pi^x$
and the walk does not visit any point in $\Pi\setminus\Pi^x$.
The opposite is also true, i.e.\
 if the walk on $\Pi$ starting from $(x,0)$ coincides with the walk on $\Pi^x$,
then $\vert W_x\vert<\infty$.
The same holds also for $s<0$: if 
$\vert W_x\vert<\infty$ and $L^0(x)<W_x$,  for  $x\in\R$, then the walk on $\Pi$ starting from $(x,0)$ 
does not visit $\Pi\setminus\Pi^x$
and if the walk does not visit $\Pi\setminus\Pi^x$ then $\vert W_x\vert<\infty$.

\begin{lemma}
 \label{xi_empty_shift}
If $s>0$ and if $\Xi$ is almost surely the empty set, then 
$\Psi=\{X\in\Pi^0:\vert W_X\vert <\infty\}$ is almost surely the empty set.
\end{lemma}
\proof
Since $W_X$ is identically distributed for all $X\in\Pi^0$, 
$\Psi$ is stationary and ergodic.
Suppose, on the contrary, that $\Psi$ is almost surely a non-empty set.
For $d>0$ let $\Psi_d=\{y\in\Pi^0:\vert W_x\vert <d\}$ and note that
$\bigcup_{d>0}\Psi_d=\Psi$.
Then there exists $d$ such that $\Psi_d$
is almost surely a non-empty set.

Let $\widetilde \Psi_d$ be the set of all $X\in \Psi_d$ which satisfy the following.
First, there are no points in $(X-d,X)\times\{0,r\}$. 
Secondly, there is $Y\in\Pi^0$, $Y<X$, such that 
the distance between $Y$ and $\max\{Z\in\Pi^0:Z<Y\}$ is greater than $\frac{r^2+s^2}{2s}$.
Thirdly, 
the walk starting from $(Y,0)$ stays in $(Y,\infty)\times\{0,r\}$ until it visits $(X,0)$
and it never visits 
$(Y+s,r)$.
Since there is a positive probability that all three conditions occur and this probability is independent of $W_x$, $\widetilde  \Psi_d$ is almost surely a 
non-empty set.
But, 
then by definition of $\Xi$, for every $X\in\widetilde \Psi_d$, $Y+s\in\Xi$
and $\Xi$ is a non-empty set, which is a contradiction.
\qed

In the next two lemmas we use the random variable $\widetilde D_X$, $X\in\Pi^r$, 
which can be compared with the corresponding random variable in Section \ref{par_indep}
defined in \eqref{tildeDy}.
For $X\in\Pi^r\setminus\Xi$ set $\widetilde D_X=0$.
For $X\in \Xi$ denote the points of $\Xi\cap(\infty,X]$ in decreasing order
$$\dots<\widetilde z_2<\widetilde z_1<\widetilde z_0=X$$
and define $\widetilde D_X$ as
\begin{align}
\label{tildeDy_shift}
\widetilde D_X =&\sup_{i\geq 0}\left\{2\widetilde z_i-\widetilde z_{i+1}-X\right\}.
\end{align}

Also we define the set $\Xi_{d_1,d_2}$ in the same way as in Section \ref{par_indep}.
For $d_1, d_2>0$ define
\begin{align*}
 \Xi_{d_1,d_2}=\{X\in\Xi:&~
\widetilde D_{X}<d_1 \text{ and there exists }Y\in\Pi^r\text{ such that }0< Y-X<d_2\text{ and }\\
&~  \widehat\Pi\cap( Y, Y+r)=\emptyset\},
\end{align*}
where $\widehat \Pi=\Pi^0\cup\Pi^r$. 

The following lemma corresponds to Lemma \ref{d_1_new}.
Since the proof is very similar, it is not included here.

\begin{lemma}
\label{d_1_new_shift}
If $s>0$ and if $\Xi$ is almost surely a non-empty set, then 
there exists $d_1,d_2>0$ such that $\Xi_{d_1,d_2}$ is  almost surely a non-empty set.
Moreover,  $\Xi_{d_1,d_2}$ is a stationary and ergodic process.
\end{lemma}

We study the greedy walk starting from the point $(0,0)$, which is almost surely not a point of $\Pi$.
From now on denote the points of $\Pi^0$ by 
 $$\ldots <X_{-2}<X_{-1}\leq 0< X_1<X_2< \dots .$$
Also, we let $\widehat \Pi=\Pi^0\cup \Pi^r$ be the shadow of all the points of the process $\Pi$.
As in Section \ref{par_indep}, we denote by $\Pi_n$ the set of points that are not visited until time $n$.
Similarly, $\Pi^0_n$, $\Pi^r_n$ denotes unvisited points of $\Pi^0$ and $\Pi^r$ until time $n$, respectively, 
and $\widehat \Pi_n=\Pi^0_n\cup\Pi^r_n$.

We define $T_A=\inf\{n\geq 0: \widehat S_n\in A\}$ to be the first time the walk visits $A\times\{0,r\}$, 
where $A$ is a subset of $\R$.

Define the variable $D_{x}$, $x\in \R$, $x>0$, 
as follows.
If $T_{(-\infty,0)}<T_{[x,\infty)}$ then $D_{x}=0$. 
Otherwise, $0< \widehat S_1,\widehat S_2,\dots,\widehat S_{T_{[x,\infty)}-1}<x,\ 
\widehat S_{T_{[x,\infty)}}\geq x$ and 
we label the remaining points of 
$\widehat \Pi_{T_{[x,\infty)}}$ in the interval $(0,x)$ by $z_1,z_2,\dots,z_{n-1}$ so that 
\begin{eqnarray*}
\label{rem_marks_shift}
 0=z_n<z_{n-1}<\dots<z_1<z_0=x.
\end{eqnarray*}
Let then
\begin{eqnarray*}
\label{D_y_shift}
 D_x=\max_{0\leq i\leq n-1}\left\{(z_i-z_{i+1})-(x-z_i)\right\}=\max_{0\leq i\leq n-1}\left\{2z_i-z_{i+1}-x\right\}.
\end{eqnarray*}
From the definition it follows that $0\leq D_x\leq x$.
As in Section \ref{par_indep}, this random variable measures how large should be the minimal distance between $(x,0)$ or $(x,r)$ and points in $(x,\infty)\times\{0,r\}$, so that the walk possibly visits a point in $(-\infty,0)\times\{0,r\}$ before visiting any point in $(x,\infty)\times\{0,r\}$.

We prove next that $D_X<d_0$ for infinitely many $X\in\Pi^0$. 
The proof is divided into two parts, one discussing the case when $\Xi$
is almost surely the empty set and another one discussing the case when $\Xi$ is a non-empty set.
The proof of the second case follows in the similar way as the second part of the proof of Lemma \ref{d_0},
so we are not going to write all the details here.

\begin{lemma}
\label{d_0_shift}
If $s>0$, then there exists $d_0<\infty$ such that, almost surely, 
$D_{X_k+s}<d_0$ and $X_{k+1}-X_k>r+s$ for infinitely many $k>0$.
\end{lemma}
\proof
Assume first that $\Xi$ is almost surely the empty set.
Then, by Lemma \ref{xi_empty_shift}, $\{X\in\Pi^0\cap (0,\infty):\vert W_X\vert <\infty\}$ is almost surely empty.
Observe that $d((X_1,0),(0,0))<d((X_{1}+s,r),(0,0))$ and if $X_{-1}+s>0$ then $d((X_{-1},0),(0,0))<s<r<d((X_{-1}+s,r),(0,0))$. 
Thus, if the walk starts from $(0,0)$ and $\widehat S_1>0$, then $S_1$ must be $(X_1,0)$.
Assume that $\{\widehat S_n>0 \text{ for all }n\geq 1\}$ occurs with positive probability.
Then one of the following three events also has positive probability.

First, $\{X_{-1}+s<0,\ \widehat S_n>0 \text{ for all }n\geq 1\}$.
But, this event implies that $\vert W_{X_1}\vert<\infty$ and $X_1\in \Psi$, which has probability $0$.

Secondly, $\{0<X_{-1}+s,\ \widehat S_n>0 \text{ for all }n\geq 1 \text{ and the walk does not visit }(X_{-1}+s,r)\}$. 
If this event occurs then the walk does not visit any point in $(0,X_{-1}+s)\times\{r\}$,
because $X_{-1}+s<r/\sqrt{3}$ and 
the distance from any point in the interval $(0,X_{-1}+s)\times\{r\}$ to the point $(X_{-1}+s,r)$ 
is smaller than distance from any point in $(0,X_{-1}+s)\times\{r\}$ to a point in $(0,\infty)\times\{0\}$.
Therefore, from $S_1=(X_1,0)$ the walk visits just the points of $\Pi^{X_1}$ and we can conclude that 
$\vert W_{X_1}\vert <\infty$, 
which is impossible.

Thirdly, $\{\widehat S_n>0 \text{ for all }n\geq 1\text{ and the walk visits }(X_{-1}+s,r)\}$.
Assume that this event occurs.
There are almost surely finitely many points in $(0,X_{-1}+s]\times\{r\}$ and thus there is 
a time $k$ when that interval is visited for the last time.
From Lemma \ref{empty_int} it follows that 
all points in $(\widehat S_k,\max_{0\leq i\leq k}\widehat S_i]\times\{0,r\}$ are visited up to time $k$.
Therefore, $S_{k+1}$ is in $(\max_{0\leq i\leq k}\widehat S_i,\infty)\times\{0,r\}$.
Since the distance from $S_k$ to $(X_{-1},0)$ is at most $\sqrt{r^2+s^2}$,
the distance from $S_k$ to $S_{k+1}$ is less than $\sqrt{r^2+s^2}$.
Thus, we can conclude that both $S_k$ and $S_{k+1}$ belong to the cluster around $(X_1,0)$.
Moreover, the greedy walk did not visit another cluster before time $k$ and it moved from 
line $0$ to line $r$ only once.

At time $k$ there might be some unvisited points of the cluster around $(X_1,0)$ on line $r$
whose shifted copies are visited before time $k$.
 Those points are visited directly after $S_k$, because those points are closer to $S_k$
than any unvisited point on line $0$.
After visiting those points, all remaining unvisited points of the cluster around $(X_1,0)$ have unvisited shifted copy.
We can think about that part of the cluster as a new cluster.
The walk visits the next cluster starting from the indented or the unindented leading point.
If it visits first the indented point, then the walk visits consecutively all the points of that cluster
and afterwards it visits the unindented leading point of the next cluster.
Once the walk is at the unindented leading point $(Y,0)$, 
all points in $(0,Y)\times\{0,r\}$ are visited except possibly some points in 
$(0,X_{-1}+s]\times\{r\}$ which are never visited.
Since $\widehat S_n>0$  for all $n\geq 1$, 
we can conclude that the walk after visiting $(Y,0)$ stays in $\Pi^Y$.
But, then $\vert W_Y\vert<\infty$ and $\Psi$ is not empty, which is a contradiction.

Since these three events almost surely do not occur, 
also $\{\widehat S_n>0 \text{ for all }n\geq 1\}$ does not occur. 
Thus $T_{(-\infty,0)}<\infty$, almost surely.
This together with the definition of $D_x$, implies that $D_x=0$ for all large enough $x\geq 0$. 
Since $X_{k+1}-X_k>r+s$ for infinitely many $k>0$, the 
claim of the lemma holds for any $d_0>0$.

If $\Xi$ is a non-empty set, the proof follows in the same way as 
the corresponding part of the proof of Lemma \ref{d_0}. 
Thus we omit the proof here and we only emphasize that points $\Xi\cap(0,\infty)$ are never visited 
because the condition $X-L^0(X)>\frac{r^2+s^2}{2s}$ implies that for $X\in\Xi\cap(0,\infty)$
points in $(-\infty,X-s)\times\{0,r\}$ are closer to $(X-s,0)$ than to $(X,r)$. 
Thus, the point $( X-s,0)$ is visited first and then by the definition of $\Xi$ the walk never visits $(X,r)$.
\qed

Using that event $D_{X_k+s}<d_0$ occurs for infinitely many $k>0$ and 
$X_{k+1}-X_k>d_0$ occurs for infinitely many $k>0$, we show in the next lemma that 
there are infinitely many $k>0$ such that both events occur simultaneously. 

\begin{lemma}
\label{prob_is_1_shifted}
If $s>0$ then, almost surely, the events 
$$A_k=\{X_{k+1}-X_k>D_{X_k+s}-X_{-1}+r+s\}$$
 occur for infinitely many $k>0$.
\end{lemma}
\proof
Let $j_{s}=\max\{i: X_i<-s\}$, $B_k=\{X_{j_s},X_{j_s+1},\dots,X_{-1},X_1,\dots,X_{k-1},X_k\}$ and
$\mathcal{F}_k=\sigma(B_k)$.
Let  $T_A^{\sigma}$ and $D_{x}^{\sigma}$ be the analogues of  $T_A$ and $D_{x}$
for the greedy walk on the set of points $\left(B_k\times\{0\}\right)\cup\left((B_k+s)\times\{r\}\right)$.
When $T_{[X_k+s,\infty)}<T_{(-\infty,0)}$,
the greedy walk on $\Pi$ and the walk on the restricted set are the same until time $T_{[X_k+s,\infty)}$.
Moreover, if $X_{k+1}-X_k>r+s$ then $\widehat S_{T_{[X_k+s,\infty)}}=X_k+s$, $T_{X_k+s}^{\sigma}= T_{[X_k+s,\infty)}$ and
 $D_{X_k+s}^{\sigma}= D_{X_k+s}$.

Let $A_k^{\sigma}=\{X_{k+1}-X_k>D_{X_k+s}^{\sigma}-X_{-1}+r+s\}$
and observe that $A_k^{\sigma}\in \mathcal F_{k+1}$.
For $d_0>0$ we have 
\begin{eqnarray*}
 \P\left(A_k^{\sigma} ~\vert~\mathcal F_k\right)&\geq& \P\left(D_{X_k+s}^{\sigma}<d_0,X_{k+1}-X_k>d_0-X_{-1}+r+s ~\vert~\mathcal F_k\right)\\
&=&\1_{\{D_{X_k+s}^{\sigma}<d_0\}}\P\left(X_{k+1}-X_k>d_0-X_{-1}+r+s ~\vert~\mathcal F_k\right)\\
&=& \1_{\{D_{X_k+s}^{\sigma}<d_0\}}e^{-(d_0-X_{-1}+r+s)} \quad \text{a.s.}
\end{eqnarray*}
The first equality above holds because $\{D_{X_k+s}^{\sigma}<d_0\}\in \mathcal F_k$.
The second equality follows from the facts that $X_{-1}\in \mathcal F_k$ and $X_{k+1}-X_k$ is exponentially distributed with mean 1
and independent of $\mathcal F_k$.
By Lemma \ref{d_0_shift}, we can choose $d_0$ such that $D_{X_k+s}<d_0$ and $X_{k+1}-X_k>r+s$ 
for infinitely many $k$, almost surely. Since, $D_{X_k+s}^{\sigma}= D_{X_k+s}$ whenever $X_{k+1}-X_k>r+s$, 
also $D^{\sigma}_{X_k+s}<d_0$ for infinitely many $k$ and, thus, 
$$\sum_{k=1}^{\infty}\P\left(A_k^{\sigma}~\vert~\mathcal F_k\right)=\infty\quad \text{a.s.}$$ 
It follows now from Lemma \ref{extended} that 
$$\P(A_k^{\sigma} \text{ for infinitely many } k\geq 1)=1.$$
Since $A_k^{\sigma}\subset \{X_{k+1}-X_k>r+s\}$ and $A_k=A_k^{\sigma}$ whenever $X_{k+1}-X_k>r+s$, also
$$\P(A_k \text{ for infinitely many } k\geq 1)=1.$$
\vskip-7mm\hfill$\Box$\vskip3mm

Whenever $A_k$ occurs, the greedy walk is forced to visit $(-\infty,0)\times\{0,r\}$
before visiting $[X_{k+1},\infty)\times\{0,r\}$. 
This together with Lemma \ref{prob_is_1_shifted} immplies  that $T_{(-\infty,0)}<\infty$ when $s>0$.
The same is also true if $s<0$ and to prove this we use that $T_{(-\infty,0)}<\infty$, almost surely, for $s>0$.

\begin{lemma}
\label{povratak_shifted}
Almost surely, $T_{(-\infty,0)}<\infty$.
\end{lemma}
\proof
When $s>0$ 
one can show in the same way as in Lemma \ref{povratak}
that whenever $A_k=\{X_{k+1}-X_k>D_{X_k+s}-X_{-1}+r+s\}$ occurs, the walk visits $(-\infty,0)\times\{0,r\}$ before visiting $[X_{k+1},\infty)\times\{0,r\}$.
By Lemma \ref{prob_is_1_shifted}, the event $A_k$ occurs for some $k$, almost surely, and hence
$T_{(-\infty,0)}<\infty$, almost surely.

Furthermore,
\begin{equation}
\label{condition3}
 \P(T_{(-\infty,0)}<\infty ~\vert ~X_{-1},X_1)=1 \quad\text{a.s.},
\end{equation}
for any absolutely continuous distribution of $(X_{-1},X_1)$ on $(-\infty,0)\times(0,\infty)$
which is independent of
$\Pi^0\cap(X_1,\infty)$.

Assume now on the contrary that  $\P(T_{(-\infty,0)}=\infty)>0$ for $s<0$.
If $T_{(-\infty,0)}=\infty$, 
then $S_1$ is in $(0,\infty)\times\{0\}$
or in $(0,\infty)\times\{r\}$.
If $S_1$ is on line $0$ then $S_1$ is indented and
the walk consecutively visits all points of its cluster.
Since the walk stays in 
$(0,\infty)\times\{0,r\}$, the last visited point of this cluster is $(X_1+s,r)$ and $X_1+s>0$.
Let $T$ be the time when the walk visits $(X_1+s,r)$ and let $Y_0=(X_1+s,r)$.
Moreover, let $Y_{-1}=X_{-1}+s$ and
let $Y_1=\min\{Y\in\Pi^r:Y>\max_{0\leq n\leq T} \widehat S_n\}$ be the leading indented point of the next cluster on the right of $Y_0$.
If $S_1$ is on line $r$, then the greedy walk is the same as if the walk starts from $(0,r)$.
Thus let $Y_0=(0,r)$, $T=0$, $Y_{-1}=X_{-1}+s$ and $Y_{1}=X_{1}+s$.

From \eqref{condition3} it follows that the walk on $\theta_{Y_0}(\Pi_T)$ 
starting from $(0,0)$ and given $X_{-1}=Y_{-1}$ and $X_1=Y_1$, visits $(-\infty,0)\times\{0,r\}$ in a finite time, almost surely.
In other words, the walk on 
$\Pi_T$ starting from $Y_0$ visits $(X_{-1},0)$ or $(X_{-1}+s,r)$ in a finite time, almost surely.
This contradicts the assumption that $\P(T_{(-\infty,0)}=\infty)>0$.
Therefore, the claim of the lemma holds also for $s<0$.
\qed

\proof[Proof of Theorem \ref{thm_emptyset4}]
From Lemma \ref{povratak_shifted} it follows that for any $\vert s\vert<r/2$
\begin{equation}
\label{condition2}
 \P(T_{(-\infty,0)}<\infty ~\vert ~X_{-1},X_1)=1 \quad\text{a.s.},
\end{equation}
for any absolutely continuous distribution of $(X_{-1},X_1)$ on $(-\infty,0)\times(0,\infty)$
which is independent of
$\Pi^0\cap(X_1,\infty)$.

We prove the theorem for $s>0$. The proof of the theorem for $s<0$ follows in a similar way. 
Let us first look at the cluster around the starting point of the greedy walk $(0,0)$.
If $X_1-X_{-1}>\sqrt{r^2+s^2}$ this cluster is empty.
If the cluster is not empty, it has finitely many points and 
the walk visits a point in another cluster  in a finite time.
Let $T_0$ be the first time the walk visits a point in another cluster 
($T_0=1$ if the cluster around $(0,0)$ is empty).

We assume at the moment that $\widehat S_{T_0}>0$.
For $i\geq 1$ let  
$T_i$ be the first time the greedy walk
visits $(-\infty,\min_{0\leq n\leq T_{i-1}}\widehat S_n)\times\{0,r\}$ for $i$ odd
and $(\max_{0\leq n\leq T_{i-1}}\widehat S_n,\infty)\times\{0,r\}$ for $i$ even.
That is, for $i$ odd (even) $T_i$ is the time when the walk visits 
the part of $\Pi$ on the left (right) of the vertical line $\{0\}\times\R$
 which is not visited up to time $T_{i-1}$.

Let $Y_0=S_{T_0-1}$ be the last visited point in the cluster around $(0,0)$ 
before the first visit to another cluster and let $\iota$ be the line of $S_{T_0-1}$. 
Moreover, let $Y_{-1}$ and $Y_{1}$ be the closest not yet visited points of $\Pi^{\iota}$ to $Y_0$, 
such that their shifted copy is also not visited, that is
 $$Y_{-1}=\max\{Y\in\Pi^{\iota}_{T_0}:Y<\min_{0\leq n< T_0} \widehat S_n,Y+(-1)^{{\text{\large $\mathfrak 1$}}_{r}(\iota)}\cdot s\in \Pi^{r-\iota}_{T_0}\}$$ and 
$$Y_1=\min\{Y\in\Pi^{\iota}_{T_0}:Y>\max_{0\leq n< T_0} \widehat S_n,Y+(-1)^{{\text{\large $\mathfrak 1$}}_{r}(\iota)}\cdot s\in \Pi^{r-\iota}_{T_0}\}.$$
Let 
$$I_1=\big((-\infty,Y_{-1})\times\{\iota\}\big)\cup \big((-\infty,Y_{-1}+(-1)^{{\text{\large $\mathfrak 1$}}_{r}(\iota)}\cdot s)\times\{r-\iota\}\big)$$
and 
$$I_2=\big((Y_{1},\infty)\times\{\iota\}\big)\cup \big((Y_1+(-1)^{{\text{\large $\mathfrak 1$}}_{r}(\iota)}\cdot s)\times\{r-\iota\}\big).$$
Because of the strong Markov property the distribution of  
$\Pi$ in $I_1$ and $I_2$ 
is independent of the points of $\Pi$ outside these sets.
Now let $\Pi'=\Pi\cap({I_1}\cup I_2)$.
From \eqref{condition2}, we know that the walk on $\sigma\theta_{Y_0}(\Pi')$ 
starting from $(0,0)$ and given $X_{-1}=Y_{-1}$ and $X_1=Y_1$ visits $(-\infty,0)\times\{0,r\}$ in a finite time, almost surely.
Hence, the walk on $\Pi_{T_0}$ starting at $S_{T_0}$ visits
$I_1$ or a point of the cluster around $(0,0)$ in almost surely finite time.
Denote that time by $T_1'$.

If $S_{T_1'}$ is in $I_1$, then $T_1=T_1'$.
Otherwise, $S_{T_1'}$ is in $\left((-\infty,Y_0)\times\{0,r\}\right)\setminus I_1$
and from the definition of $I_1$ 
we can deduce that the shifted copy of the point $S_{T_1'}$ must have been visited before $T_1$.  
Set now $Y_0=S_{T_1'}$ and redefine $Y_{-1}$, $Y_1$, $I_2$ and $\Pi'$ with respect to the time $T_1'$ 
instead of $T_0$.

Observe that, by Lemma \ref{empty_int}, at time $T_1'$ the set
$(\widehat S_{T_1'},\max_{0\leq n\leq T_1'} \widehat S_n)\times\{0,r\}$ is empty.
Moreover, if there are some points in $(\max_{0\leq n\leq T_1'} \widehat S_n,\infty)\times\{0,r\}$
whose shifted copy is visited, then these points are on  line $r$ and belong to one cluster.
The closest point with a still unvisited shifted copy is at a horizontal distance of at least $r-s$ from those points. 
Since the distance from $S_{T_1'}$ to $(\max_{0\leq n\leq T_1'} \widehat S_n,\infty)\times\{0,r\}$ is at least 
$\sqrt{r^2+s^2}$, the remaining points on line $r$ which do not have a shifted copy 
are closer to $S_{T_1'}$ than the closest point on line $0$.
Thus, these points, whose shifted copies are already visited, are visited before visiting $(Y_1',\iota)$ or its shifted copy.

Now, we can conclude that the walk starting from $S_{T_1'}$ visits in almost surely finite time
$I_1$ or a point of $\Pi\setminus (I_1\cup I_2)$, that is one of the remaining points of the cluster around $(0,0)$ 
or a point in  $(\max_{0\leq n\leq T_1'} \widehat S_n,Y_1')\times\{r\}$.
There are finitely many points in $\Pi\setminus (I_1\cup I_2)$ and every time the walk visits one of these points, 
we redefine $Y_0$, $Y_{-1}$, $Y_1$, $I_2$ and $\Pi'$,
and repeat the same arguments as above.
Thus the walk visits $I_1$ in almost surely finite time.

Assume now that $T_i$ is finite for some odd $i\geq 1$.
Let $Y_0=\widehat S_{T_i-1}$ and define $Y_{-1}$, $Y_1$, $I_1$, $I_2$ and $\Pi'$ as before.
Then points of $\Pi_{T_i}\setminus\Pi'$ are in $(Y_{-1},0)\times \{r\}$, $(0,Y_1)\times\{0\}$
or the cluster around $(0,0)$ and there are  almost surely finitely many such points.
By the observation above, the greedy walk visits  all points in $(Y_{-1}, Y_0)\times \{0\}$ before 
it visits $I_1$ and it visits all points in $(Y_0,Y_1)\times\{r\}$ before it visits a point in $I_2$.
Denote the time when the walk visits a point of $\Pi_{T_i}\setminus\Pi'$ before visiting $I_2$ with $T_i'$.
Set then $Y_0=S_{T_i'}$ and redefine $Y_{-1}'$, $Y_1$ and $\Pi'$, with the respect with the time $T_i'$.
Again, by \eqref{condition2}, the walk on $\sigma\theta_{Y_0}\Pi'$ 
visits $(-\infty,0)\times\{0,r\}$ in almost surely finite time. 
Thus, the walk on $\Pi_{T_i'}$ in almost surely finite time visits $I_2$ or another point in $\Pi_{T_i}\setminus\Pi'$.
Repeating this arguments for every visited point in $\Pi_{T_i}\setminus\Pi'$, 
we can see that the walk eventually visits $I_2$ and that $T_{i+1}$ is almost surely finite.

Similarly, one can show that if $T_i$, $i\geq 2$ even, is finite, 
then $T_{i+1}$ is also almost surely finite.
Therefore, inductively we can conclude that the walk almost surely crosses the vertical line $\{0\}\times\R$ infinitely many times
and, thus, it eventually visits all points of $\Pi$.
%
\qed

\begin{remark}
 We conjecture that Theorem \ref{thm_emptyset4} holds also for 
$\vert s\vert >r/2$.
For those $s$ the idea to cluster the points of $\Pi$ does not work in the same way. 
For example, the greedy walk does not always visit all points of the cluster when
it starts from the leading indented point of the cluster.
Thus, the walk more often does not visit all points of a cluster successively and 
we expect that the points that are not visited during the first visit of a cluster 
cause the walk to return and to cross the vertical line $\{0\}\times\R$ infinitely often.
\end{remark}

\vspace{0.1cm}
\acknowledgments
The author thanks Svante Janson, Takis Konstantopoulos and Erik Th\"ornblad for valuable comments.

\small
\bibliography{two_lines}
\bibliographystyle{apt}

\vspace*{0.2cm}
\begin{minipage}[]{13.5cm}
\footnotesize \sc
Katja Gabrysch, 
Department of Mathematics,
Uppsala University,
PO~Box~480, 751~06~Uppsala,
Sweden
\\
{\em E-mail address:} {\tt \href{mailto:katja@math.uu.se}{katja@math.uu.se}}
\end{minipage}

\end{document}